\documentclass[preprint,sort&compress,final,10pt]{elsarticle}

\usepackage{amsmath}
\usepackage{amssymb}
\usepackage{graphicx}
\usepackage{latexsym}
\usepackage{epstopdf}
\usepackage{natbib}
\usepackage{color}
\usepackage{hyperref}

\newtheorem{theorem}{Theorem}[section]

\newtheorem{definition}[theorem]{Definition}
\newtheorem{lemma}[theorem]{Lemma}

\newtheorem{remark}[theorem]{Remark}
\numberwithin{equation}{section}

\def\Proof{\noindent{\bf Proof.}~}
\def\qed{\hfill$\square$\smallskip}

\journal{\empty}
\date{}

\begin{document}

\begin{frontmatter}

\title{Entire Solutions for the Classical Competitive Lotka-Volterra System with Diffusion in the Weak Competition Case}

\author[au1,au2]{Yang Wang}

\address[au1]{ School of Mathematical Sciences, Shanxi University, Shanxi 030006, P.R. China.}

\ead[au1]{ywang2005@sxu.edu.cn}

\author[au1]{Xiong Li\footnote{ Partially supported by the NSFC (11571041) and the Fundamental Research Funds for the Central Universities. Corresponding author.}}

\address[au2]{School of Mathematical Sciences, Beijing Normal University, Beijing 100875, P.R. China.}

\ead[au2]{xli@bnu.edu.cn}

\begin{abstract}
In this paper we are concerned with the entire solutions for the classical competitive Lotka-Volterra system with diffusion in the weak competition. For this purpose we firstly analyze the asymptotic behavior of traveling front solutions for this system connecting the origin and the positive equilibrium. Then, by using two different ways to construct pairs of coupled super-sub solutions of this system, we obtain two different kinds of entire solutions. The construction of the first kind of entire solutions is based on these fronts, and their reflects as well as the solutions of the system without diffusion. One component of the solution starts from 0 at $t\approx-\infty$, and as $t$ goes to $+\infty$, the two component of the entire solution will eventually stay in a conformed region. Another kind of entire solutions is related to some traveling front solutions of scalar equations.
\end{abstract}


\end{frontmatter}

\section{Introduction}
In this paper, we are concerned with the classical competitive Lotka-Volterra system with diffusion
\begin{equation}\label{eq:LV}
\left\{\begin{array}{ll}
\partial_tu=\partial_{xx}u+(1-u-k_1v)u,\\[0.2cm]
\partial_tv=d\,\partial_{xx}v+r(1-v-k_2u)v,
\end{array}
\right.x\in\mathbb{R}
\end{equation}
where $k_1$, $k_2$, $r$, $d$ are positive constants and $u(x,t)$, $v(x,t)$ denote the population density of two competitive species that are nonnegative.

To begin with this paper, we remark that, as stated in \cite{mt09}, the solutions $(u(t),v(t))$ of \eqref{eq:LV} without diffusion exhibit the following asymptotic behavior as $t\rightarrow+\infty$:\\
(i) \ if $0<k_1<1<k_2$, then $(u(t),v(t))\rightarrow(1,0)$ ($u$ survives);\\
(ii)\ if $0<k_2<1<k_1$, then $(u(t),v(t))\rightarrow(0,1)$ ($v$ survives);\\
(iii)\ if $k_1$, $k_2>1$, then $(u(t),v(t))\rightarrow(1,0)$ or $(u(t),v(t))\rightarrow(0,1)$ depending on the initial condition (strong competition and bistability);\\
(iv)\ if $0<k_1,k_2<1$, then $(u(t),v(t))$ converges to the positive equilibrium (weak competition, $u$ and $v$ coexist).

Furthermore these results can be extended to system \eqref{eq:LV} (\cite{m79}). In this paper, we only pay attention to the existence of entire solutions and other properties for system \eqref{eq:LV} under the case (iv): $0<k_1,k_2<1$, namely the weak competition case, since there are relatively abundant results under the cases (i)-(iii), which will be depicted at length in the following. In this case, the above system has four equilibria that are $(0,0)$, $(1,0)$, $(0,1)$ and $(u^{\ast},v^{\ast}):=\left(\frac{1-k_1}{1-k_1k_2},
\frac{1-k_2}{1-k_1k_2}\right)$. Moreover, we also remark that $u^{\ast}+k_1v^{\ast}=k_2u^{\ast}+v^{\ast}=1$, which will be used in the sequel.

A great deal of papers focus on the study of traveling wave solutions of system \eqref{eq:LV} such as \cite{al91,alt08,cg84,lhl08,lhf11,tf80,v95,yy11} and references therein, which is a significant and of particular interesting issue in reaction diffusion systems. Moreover, many authors paid more attention to the existence of monotone traveling wave solutions, namely traveling front solutions.

In \cite{tf80}, the authors discussed the existence of traveling front solutions and traveling wave solutions for a general Lotka-Volterra competition model
\begin{equation}\label{eq:GLV}
\left\{\begin{array}{ll}
\partial_tu=\partial_{xx}u+uf(u,v),\\[0.2cm]
\partial_tv=d\,\partial_{xx}v+vg(u,v),
\end{array}
\right.x\in\mathbb{R}
\end{equation}
where the functions $f$ and $g$ satisfy the following assumptions:\\
(a)  $f$, $g\in C^1$ have a positive zero $(u^{\star},v^{\star})$, that is, $f(u^{\star},v^{\star})=g(u^{\star},v^{\star})=0$;\\
(b) if $0<u<u^{\star}$, $0<v<v^{\star}$, then $0<f(u,v)<f(0,0)$, \ $0<g(u,v)<g(0,0)$;\\
(c) if $0<u<u^{\star}$, $0<v<v^{\star}$, then $\partial_{u}f(u,v)<0$, $\partial_{v}f(u,v)\leqslant0$, $\partial_{u}g(u,v)\leqslant0$, $\partial_{v}g(u,v)<0$;\\
(d) the eigenvalues of the matrix
$$
\begin{pmatrix}u^{\star}\partial_{u}f(u^{\star},v^{\star}) & u^{\star}\partial_{v}f(u^{\star},v^{\star})\\[0.3cm] v^{\star}\partial_{u}g(u^{\star},v^{\star}) & v^{\star}\partial_{v}g(u^{\star},v^{\star})
\end{pmatrix}
$$
have negative real parts. As stated in \cite{al91}, Perron-Frobenius theorem and the assumption (c) imply that the matrix in assumption (d) has a real and negative eigenvalue. Hence, both eigenvalues are real. Furthermore in the papers \cite{al91} and \cite{alt08}, the authors extended the results in \cite{tf80} into $N$-equations. Moreover, there exists a family of planar front solutions for \eqref{eq:GLV} on $\mathbb{R}^{n}$ with different types of reaction terms (\cite{v95}).

In fact, it is easy to verify that the reaction terms in \eqref{eq:LV} under the case (iv) fully satisfy the above assumptions. Therefore, from \cite{tf80}, for $c\geqslant c_{\min}:=2\max\{1,\sqrt{rd}\}$, \eqref{eq:LV} admits a family of traveling front solutions connecting $(0,0)$ and $(u^{\ast},v^{\ast})$. More precisely, the traveling front solution $(u(x,t),v(x,t))=(\phi(\xi),\psi(\xi))$ ($\xi=x+ct$, $c\geqslant c_{\min}$) for \eqref{eq:LV} satisfies
\begin{equation}\label{eq:OLV}
\left\{\begin{array}{ll}
\phi''-c\phi'+(1-\phi-k_1\psi)\phi=0,\\[0.2cm]
d\psi''-c\psi'+r(1-\psi-k_2\phi)\psi=0
\end{array}
\right.
\end{equation}
with
\begin{align}
&\label{eq:BC}\lim\limits_{\xi\rightarrow-\infty}(\phi(\xi),\psi(\xi))=(0,0),\ \ \ \lim\limits_{\xi\rightarrow+\infty}(\phi(\xi),\psi(\xi))=(u^{\ast},v^{\ast}),\\[0.2cm]
&\label{eq:MC}\phi(\xi),\ \ \psi(\xi)>0,\hspace{2cm} \phi'(\xi),\ \ \psi'(\xi)>0.
\end{align}
We also remark that if $(\phi(x+ct),\psi(x+ct))$ is a traveling front solution of \eqref{eq:OLV}-\eqref{eq:MC}, then the reflect $(\tilde{\phi}(\tilde{\xi}),\tilde{\psi}(\tilde{\xi}))=
(\phi(-x+ct),\psi(-x+ct))\ (\tilde{\xi}=-x+ct)$ is also a traveling front solution with the opposite speed satisfying
$$\lim\limits_{\tilde{\xi}\rightarrow-\infty}
(\tilde{\phi}(\tilde{\xi}),\tilde{\psi}(\tilde{\xi}))=(u^{\ast},v^{\ast}),\ \ \ \lim\limits_{\tilde{\xi}\rightarrow+\infty}
(\tilde{\phi}(\tilde{\xi}),\tilde{\psi}(\tilde{\xi}))=(0,0),$$
and
$$
\tilde{\phi}(\tilde{\xi}),\ \ \tilde{\psi}(\tilde{\xi})>0,\hspace{2cm} \tilde{\phi}'(\tilde{\xi}),\ \ \tilde{\psi}'(\tilde{\xi})<0.
$$
That is, if \eqref{eq:OLV}-\eqref{eq:MC} admits a traveling front solution, then a traveling front solution exists with the opposite speed simultaneously.

However, it is not enough to understand the dynamical structure of solutions of \eqref{eq:LV} by only considering traveling wave solutions. Recently, the existence of entire solutions, which are classical solutions and defined for all $(x,t)\in\mathbb{R}\times\mathbb{R}$, has been widely discussed.

In \cite{hn99}, Hamel and Nadirashvili dealt with KPP equation
\begin{align}
&\label{eq:KPP-Fisher}\partial_tu=\partial_{xx}u+f(u),\ \ (x,t)\in\mathbb{R}\times\mathbb{R},\\[0.2cm]
&\label{eq:MOS}f'(0)>0\ \ \ \ \ \ f'(1)<0.
\end{align}
The existence of entire solutions is proved by the comparison theorem and super-sub estimates, which consists of traveling front solutions and solutions to the diffusion-free system. Moreover, they also pointed out that the solutions to \eqref{eq:KPP-Fisher} depending only on $t$ and traveling wave solutions are typical examples of entire solutions and showed various entire solutions of \eqref{eq:KPP-Fisher} with \eqref{eq:MOS} in their subsequent paper \cite{hn01}. More importantly, from the geometrical point of view, in \cite{hn99}, they indicated that \eqref{eq:KPP-Fisher} has a 2-dimensional manifold of entire solutions of traveling wave type, which are $\tilde{\phi}(x+\tilde{c}t+\tilde{h})$ and $\tilde{\phi}(-x+\tilde{c}t+\tilde{h})$, where $\tilde{h}$ varies in $\mathbb{R}$ and $\tilde{c}$ varies in $[c_{\ast},+\infty]$ with $c_{\ast}=2\sqrt{f'(0)}$. They also established 5-dimensional, 4-dimensional, 3-dimensional manifolds of entire solutions and showed that each 2-dimensional manifold of entire solutions of traveling wave type is on the boundary of a 3-dimensional manifold of entire solutions, which is also on the boundary point of a 5-dimensional one. As stated in \cite{cg05}, after a space and time translation, these manifolds can be reduced to 3, 2, 1-dimensional manifolds, where the parameters vary in $[c_{\ast},+\infty]\times[c_{\ast},+\infty]\times[-\infty,+\infty]$.

While for the bistable case, namely both $f'(0)<0$ and $f'(1)<0$, Yagisita in \cite{y03} revealed that the annihilation process is approximated by a backward global solution of \eqref{eq:KPP-Fisher}, which is an entire solution. For Allen-Cahn equation
$$\partial_tu=\partial_{xx}u+u(1-u)(u-a)$$
with $a\in(0,1)$, which is a special example of \eqref{eq:KPP-Fisher} in \cite{y03}, Fukao, Morita and Ninomiya in \cite{fmn04} proposed a simple proof for the existence of entire solutions, which was already found in \cite{y03}, by using the super-sub solution method and the exact traveling front solutions. Moreover, Guo and Morita in \cite{gm05} extended the conclusions in \cite{hn99} and \cite{y03} to more general case. In addition, Chen and Guo in \cite{cg05} used a quite different method to construct the super-sub solutions to obtain the similar results for this more general system in \cite{gm05}.

From the dynamical view the study of entire solutions is essential for a full understanding of the transient dynamics and the structures of the global attractor as mentioned in \cite{mn06}. Recently, there has been large numbers of papers about the existence of entire solutions of scalar equation, for example, see \cite{llw08,lww08,wlr09,ct12,wlw09} and the references therein.

In 2009, Morita and Tachibana in \cite{mt09} firstly extended the existence of entire solutions from scalar equations into system \eqref{eq:LV} under the cases (i), (ii) and (iii) by employing similar ideas in \cite{cg05}, \cite{fmn04}, \cite{gm05} and \cite{mn06}. This entire solution behaves as two traveling front solutions coming from both sides of the $x$-axis at $t\approx-\infty$ and one component converges to $1$ while the other converges to $0$, as $t$ goes to $+\infty$. With similar methods, in \cite{wl13} Wang and Lv obtained entire solutions of system \eqref{eq:LV} by changing \eqref{eq:LV} into the cooperative system with the extra condition $d=1$ and (iv), and also obtained entire solutions of the classical Lotka-Volterra cooperative system. The asymptotic behavior of entire solutions in \cite{wl13} is similar as that in \cite{mt09}. In addition, in \cite{mt09} the existence of entire solutions of \eqref{eq:LV} is based on traveling front solutions which connect $(0,1)$ and $(1,0)$, while in \cite{wl13} the existence of entire solutions of \eqref{eq:LV} depends on traveling front solutions which connect $(0,1)$ and the positive equilibrium. In addition, entire solutions of the cooperative system in \cite{wl13} are related to traveling front solutions connecting the origin and the positive equilibrium.

In addition, there is a technical condition
$$\frac{\phi(\xi)}{1-\psi(\xi)}\geqslant\theta_0>0$$
for the existence of entire solutions in \cite{mt09}. In \cite{wl15}, Wang and Li showed some sufficient and necessary conditions for this technical condition and partially proved the result still holds without this condition. In addition, except for the above mentioned papers, for the existence of entire solutions of Lotka-Volterra system one can see \cite{lzz15}, \cite{l12} and \cite{wl10} for more details.

From the above statement, we remark that there is no any results for the existence of entire solutions for system \eqref{eq:LV} in the case (iv) based on the solutions to \eqref{eq:OLV} with \eqref{eq:BC} and \eqref{eq:MC}. In this paper, encouraged by \cite{mt09} and \cite{wl13}, we initially want to find entire solutions which can demonstrate that one species invades from both sides of $x-$axis and coexists with the other at last. However, a pair of coupled super-sub solutions constructed in \cite{wl13} for the cooperative system can not be directly used in the competitive system, due to the different monotonicity of these systems. Thus we firstly use the traveling front solutions of the above system connecting the origin and the positive equilibrium and their reflects as well as the solutions of the above system without diffusion to construct different pairs of coupled  super-lower solutions leading to the existence of different kinds of entire solutions. One of them has the following asymptotic behavior. One component of this entire solution start from 0 at $t\approx-\infty$. As $t$ converges to $+\infty$, two component of the entire solution will ultimately stay in a conformed region. This phenomenon implies one species invades from both sides of $x-$axis and will mix with the other. In addition, this entire solution exhibits quite different behavior compared with that in \cite{mt09} and \cite{wl13}.

For the sake of realizing the initial conjecture and finding more types of entire solutions, different methods have been employed. With the idea coming from \cite{lhl08,lhf11}, we construct the pair of coupled super-sub solutions related to traveling front solutions of some scalar equations, then another kind of entire solutions is found.

Moreover, from their applications, we note that this method depends on the asymptotic behavior of traveling front solutions connecting the origin and the positive equilibrium as $\xi\rightarrow\pm\infty$. In addition, in the proof of the existence of entire solutions for \eqref{eq:LV}, we use the super-sub solution method and then need to estimate the asymptotic behavior and the boundness of the constructed pair of coupled super-sub solutions, which are related to the asymptotic behavior of traveling front solutions connecting the origin and the positive equilibrium as $\xi\rightarrow\pm\infty$ as well. Therefore, it is significant for us to study the asymptotic behavior of traveling front solutions for \eqref{eq:OLV} at $(0,0)$ and $(u^{\ast},v^{\ast})$, respectively, which will be achieved by linearizing \eqref{eq:OLV} and the stable and unstable manifold theorem.

Though there are some results about the asymptotic behavior of the above traveling front solutions, we employ the above method to obtain full and accurate conclusions. For instance, with the aid of Laplace transform, the asymptotic behavior of traveling front solutions for $c>c_{\min}$ was established in the paper \cite{l08}, which is
$$\lim\limits_{\xi\rightarrow-\infty}\Bigl(\phi(\xi)e^{-\tilde{\lambda}_1(\xi+\tilde{h}_1)},\,
\psi(\xi)e^{-\tilde{\lambda}_2(\xi+\tilde{h}_2)}\Bigr)=(1,1)$$
for some $\tilde{\lambda}_i$, $\tilde{h}_i>0$, $i=1,2$. Our results improve this conclusions and extend to the case $c=c_{\min}$ and include the asymptotic behavior as $\xi\rightarrow+\infty$. We also note that in \cite{wl13}, the authors stated the asymptotic behavior of traveling front solutions for the cooperative system as $\xi\rightarrow\pm\infty$, while, here we show more details for the competitive system and only the smaller negative eigenvalue of the linearization matrix at $(u^{\ast},v^{\ast})$ plays role in the asymptotic behavior as $\xi\rightarrow+\infty$.

The paper is organized as follows. We are devoted to deeply analyzing the asymptotic behavior of traveling front solutions connecting the origin and the positive equilibrium for \eqref{eq:LV} as $\xi\rightarrow\pm\infty$ in Section 2. The definition of a pair of coupled super-sub solutions as well as the existence and qualitative properties of entire solutions of general quasi-monotone decreasing reaction-diffusion systems are restated in Section 3. In the end, two different kinds of entire solutions for system \eqref{eq:LV} and their qualitative properties are discussed in Section 4.


\section{Asymptotic behavior of traveling front solutions}

As emphasized in the introduction, the asymptotic behavior of traveling front solutions connecting the origin and the positive equilibrium for \eqref{eq:LV} as $\xi\rightarrow\pm\infty$ should be analyzed at first. In the sequel, we always assume $(\phi(\xi),\psi(\xi))$ is a pair of traveling front solutions connecting the origin and the positive equilibrium for \eqref{eq:LV}, that is, which is the solution of \eqref{eq:OLV} with \eqref{eq:BC} and \eqref{eq:MC}. Set $\phi'=Y$ and $\psi'=Z$, then \eqref{eq:OLV} is equivalent to
\begin{equation}\label{eq:OLVS}
\left\{\begin{array}{llll}
\phi'=Y,\\[0.2cm]
Y'=cY-\phi+\phi^2+k_1\phi\psi,\\[0.2cm]
\psi'=Z,\\[0.2cm]
Z'=\frac{c}{d}Z-\frac{r}{d}\psi+\frac{r}{d}\psi^2+\frac{r}{d}k_2\phi\psi.\\[0.2cm]
\end{array}
\right.
\end{equation}

We firstly discuss the asymptotic behavior of the solution of \eqref{eq:OLV} with \eqref{eq:BC} and \eqref{eq:MC} as $\xi\rightarrow+\infty$. Linearizing \eqref{eq:OLVS} at the point $(u^{\ast},0,v^{\ast},0)$ yields that
\begin{equation}\label{eq:OLVS1}\begin{pmatrix} \phi_+' \\ Y_+' \\ \psi_+' \\ Z_+' \end{pmatrix}=A_1\begin{pmatrix} \phi_+ \\ Y_+ \\ \psi_+ \\ Z_+ \end{pmatrix},\end{equation}
where $$A_1=\begin{pmatrix} 0 & 1 & 0 & 0 \\ u^{\ast} & c & k_1u^{\ast} & 0 \\ 0 & 0 & 0 & 1 \\ \frac{r}{d}k_2v^{\ast} & 0 & \frac{r}{d}v^{\ast} & \frac{c}{d} \end{pmatrix}.$$
The characteristic equation of the matrix $A_1$ is
$$\begin{array}{ll}\lambda^4-\left(c+\frac{c}{d}\right)\lambda^3+\left(\frac{c^2}{d}
-u^{\ast}-\frac{r}{d}v^{\ast}\right)\lambda^2+\left(\frac{cr}{d}v^{\ast}
+\frac{c}{d}u^{\ast}\right)\lambda+\frac{r}{d}(1-k_1k_2)u^{\ast}v^{\ast}=0.\end{array}$$
Apparently, it is not hard to calculate the eigenvalues of the matrix $A_1$, but the distributions of the eigenvalues can be determined by the method from \cite{tf80}. Although the proof is similar to that in \cite{tf80}, for the reader's convenience, we show the proof in the following.
\begin{lemma}\label{lem:eigen}
The equilibrium $(u^{\ast},0,v^{\ast},0)$ is hyperbolic, and both the stable subspace and the unstable subspace of \eqref{eq:OLVS1} at $(u^{\ast},0,v^{\ast},0)$ are two dimensional.
\end{lemma}

\Proof To begin with the proof, we introduce a matrix
$$\Lambda(\rho)=\begin{pmatrix} 0 & 1 & 0 & 0 \\ u^{\ast} & c & \rho k_1u^{\ast} & 0 \\ 0 & 0 & 0 & 1 \\ \frac{r}{d}k_2v^{\ast} & 0 & \frac{r}{d}v^{\ast} & \frac{c}{d} \end{pmatrix}$$
with the parameter $\rho\in[0,1]$, and the corresponding characteristic equation is $F_{\rho}(\lambda)=0$,
where
$$
\begin{array}{ll}
F_{\rho}(\lambda)=\lambda^4-\left(c+\frac{c}{d}\right)\lambda^3+\left(\frac{c^2}{d}-u^{\ast}
-\frac{r}{d}v^{\ast}\right)\lambda^2+
\left(\frac{rc}{d}v^{\ast}+\frac{c}{d}u^{\ast}\right)\lambda+d(\rho),
\end{array}
$$
and $$d(\rho)=\det(\Lambda(\rho))=\frac{r}{d}u^{\ast}v^{\ast}-\frac{\rho rk_1k_2}{d}u^{\ast}v^{\ast}>0,$$
since $0<k_1,k_2<1$.

It is easy to compute that the eigenvalues of the matrix $\Lambda(0)$ are
\begin{align*}
&\mu_1=\frac{c+\sqrt{c^2+4u^{\ast}}}{2},\ \ \ \ \ \ \ \mu_2=\frac{c-\sqrt{c^2+4u^{\ast}}}{2},\\[0.2cm]
&\mu_3=\frac{c+\sqrt{c^2+4drv^{\ast}}}{2d},\ \ \ \ \mu_4=\frac{c-\sqrt{c^2+4drv^{\ast}}}{2d}.
\end{align*}
Obviously, they are real, two of them are positive and the others are negative. Meanwhile, the corresponding characteristic equation is $F_0(\lambda)=0$, where
$$F_0(\lambda)=\left(\lambda^2-c\lambda-u^{\ast}\right)\left
(\lambda^2-\frac{c}{d}\lambda-\frac{r}{d}v^{\ast}\right).$$
We also note that the function $F_0(\lambda)$ is positive for some sufficiently large $|\lambda|$, negative for some positive and negative $\lambda's$, and $F_0(0)>0$.

From the definition of $F_{\rho}(\lambda)$, the graph of the function $F_{\rho}(\lambda)$ is either upwards or downwards as $\rho$ changes. More precisely, when $\rho$ increases from 0 to 1, $d(\rho)$ decreases leading to the downwards of the graph of the function $F_{\rho}(\lambda)$. Since for each $\rho$, $F_{\rho}(0)=d(\rho)>0$, the equation $F_1(\lambda)=0$, which is the characteristic equation of the matrix $A_1$, is continuous to have two different negative roots and two different positive roots. This finishes the proof.
\qed

From Lemma \ref{lem:eigen}, we assume that $\lambda_1>\lambda_2$ are two negative eigenvalues of the matrix $A_1$. Then the corresponding eigenvectors are
$$\begin{pmatrix} 1 \\ \lambda_1 \\ \tau_1 \\ \tau_1\lambda_1 \end{pmatrix},\ \ \ \ \ \ \ \begin{pmatrix} 1 \\ \lambda_2 \\ \tau_2 \\ \tau_2\lambda_2 \end{pmatrix},$$
where $$\tau_1=\frac{\lambda_1^2-c\lambda_1-u^{\ast}}{k_1u^{\ast}},\ \ \ \ \tau_2=\frac{\lambda_2^2-c\lambda_2-u^{\ast}}{k_1u^{\ast}}.$$
We notice that $\lambda_2<\mu_2<\lambda_1$ derived from the proof of Lemma \ref{lem:eigen}, which leads to $\lambda_1^2-c\lambda_1-u^{\ast}<0$ and $\lambda_2^2-c\lambda_2-u^{\ast}>0.$ Consequently, $\tau_1<0$ and $\tau_2>0$.

Therefore every solution of \eqref{eq:OLVS1}, converging to the origin as $\xi\rightarrow+\infty$, can be given by
$$\begin{pmatrix} \phi_+ \\ Y_+ \\ \psi_+ \\ Z_+ \end{pmatrix}=C_1\begin{pmatrix} 1 \\ \lambda_1 \\ \tau_1 \\ \tau_1\lambda_1 \end{pmatrix}e^{\lambda_1\xi}+C_2\begin{pmatrix} 1 \\ \lambda_2 \\ \tau_2 \\ \tau_2\lambda_2 \end{pmatrix}e^{\lambda_2\xi},$$
where $C_1$ and $C_2$ are arbitrarily given real numbers. By applying the stable manifold theorem, as $\xi\rightarrow+\infty$, there are some constants $\alpha$ and $\beta$ such that
\begin{align*}\phi(\xi)&=u^{\ast}-\alpha e^{\lambda_1\xi}-\beta e^{\lambda_2\xi}+h.o.t.,\\
\psi(\xi)&=v^{\ast}-\alpha\tau_1 e^{\lambda_1\xi}-\beta\tau_2 e^{\lambda_2\xi}+h.o.t..
\end{align*}

Now we show that $\alpha=0$. If $\alpha\neq0$, then we firstly suppose $\alpha>0$. Because of the above equations, for any $\beta$ and sufficiently large $\xi$, $\psi>v^{\ast}$, which contradicts \eqref{eq:BC} and \eqref{eq:MC}. On the other hand, if $\alpha<0$, then for any $\beta$ and sufficiently large $\xi$, $\phi>u^{\ast}$, which also is a contradiction. Thus the asymptotic behavior can be refined as
\begin{align*}\phi(\xi)&=u^{\ast}-\beta e^{\lambda_2\xi}+h.o.t.,\\
\psi(\xi)&=v^{\ast}-\beta\tau_2 e^{\lambda_2\xi}+h.o.t.,
\end{align*}
where $\beta$ is a positive constant. As a result, the following lemma is obtained.
\begin{lemma}\label{lem:asym1}
The asymptotic behavior of the traveling front solution $(\phi,\psi)$, as $\xi\rightarrow+\infty$, is
\begin{align*}
\phi(\xi)&=u^{\ast}-\beta e^{\lambda_2\xi}+h.o.t.,\\
\psi(\xi)&=v^{\ast}-\beta\tau_2 e^{\lambda_2\xi}+h.o.t.,
\end{align*}
where $\beta$ and $\tau_2$ are two positive constants.
\end{lemma}

Next, we are devoted to investigating the asymptotic behavior of the traveling front solution $(\phi,\psi)$ when $\xi\rightarrow-\infty$, which is finished by several lemmas. As $\xi\rightarrow-\infty$, compared with the paper \cite{l08}, here we will use a totally different method similar to that in \cite{mt09} to study the asymptotic behavior and show more details, for example, the asymptotic behavior of $(\phi,\psi)$ under the case $c=c_{\min}$.

Similar to the analysis as above, the linearized system at the point $(0,0,0,0)$, is
\begin{equation*}\begin{pmatrix} \phi_-' \\ Y_-' \\ \psi_-' \\ Z_-' \end{pmatrix}=A_2\begin{pmatrix} \phi_- \\ Y_- \\ \psi_- \\ Z_- \end{pmatrix},\end{equation*}
where $$A_2=\begin{pmatrix} 0 & 1 & 0 & 0 \\ -1 & c & 0 & 0 \\ 0 & 0 & 0 & 1 \\ 0 & 0 & -\frac{r}{d} & \frac{c}{d} \end{pmatrix}.$$
Set
$$B_1:=A_2-\lambda I=\begin{pmatrix} -\lambda & 1 & 0 & 0 \\ -1 & c-\lambda & 0 & 0 \\0 & 0 & -\lambda & 1 \\ 0 & 0 & -\frac{r}{d} & \frac{c}{d}-\lambda \end{pmatrix}.$$
The characteristic equation of the matrix $A_2$ is
\begin{equation}\label{eq:chara}
\det(B_1)=\det(A_2-\lambda I)=(\lambda^2-c\lambda+1)(\lambda^2-\frac{c}{d}\lambda+\frac{r}{d})=0.
\end{equation}
Then the eigenvalues are
\begin{align*}
&\lambda_3=\frac{c+\sqrt{c^2-4}}{2},\ \ \ \ \ \ \ \ \lambda_4=\frac{c-\sqrt{c^2-4}}{2},\\[0.2cm]
&\lambda_5=\frac{c+\sqrt{c^2-4rd}}{2d},\ \ \ \ \ \lambda_6=\frac{c-\sqrt{c^2-4rd}}{2d},
\end{align*}
which are real since $c\geqslant2\max\{1,\sqrt{rd}\}$, and obviously $\lambda_3\geqslant\lambda_4>0$, $\lambda_5\geqslant\lambda_6>0$.
Then according to the different multiplicities of the eigenvalues of $A_2$, we will investigate the asymptotic behavior of $(\phi,\psi)$ as $\xi\rightarrow-\infty$ in several cases.

Firstly, we consider all eigenvalues of $A_2$ are simple roots, and then obtain the following lemma.

\begin{lemma}\label{lem:asym01}
When $\lambda_3>\lambda_4$, $\lambda_5>\lambda_6$ and $\lambda_i\neq\lambda_j$ for $i=3,4$ and $j=5,6$, the traveling front solution $(\phi,\psi)$ behaves as $\xi\rightarrow-\infty$ in the following way:
\begin{align*}
\phi(\xi)&=\alpha e^{\lambda_3\xi}+\beta e^{\lambda_4\xi}+h.o.t.,\\
\psi(\xi)&=\gamma e^{\lambda_5\xi}+\sigma e^{\lambda_6\xi}+h.o.t.,
\end{align*}
where $\beta\geqslant0$, $\sigma\geqslant0$, $\alpha>0\,(\beta=0)$, $\gamma>0\,(\sigma=0)$.
\end{lemma}

\Proof For each eigenvalue $\lambda_i$, $i=3,4,5,6$, the corresponding eigenvectors are
$$r^1_1=\begin{pmatrix} 1 \\ \lambda_3 \\ 0 \\ 0 \end{pmatrix},\ \ \ r^1_2=\begin{pmatrix} 1 \\ \lambda_4 \\ 0 \\ 0 \end{pmatrix},\ \ \ r^1_3=\begin{pmatrix} 0 \\ 0 \\ 1 \\ \lambda_5 \end{pmatrix},\ \ \ r^1_4=\begin{pmatrix} 0 \\ 0 \\ 1 \\ \lambda_6 \end{pmatrix}.$$
Then there are some arbitrarily given numbers $C_1$, $C_2$, $C_3$ and $C_4$ such that
$$\begin{pmatrix} \phi_- \\ Y_- \\ \psi_- \\ Z_-  \end{pmatrix}=C_1\begin{pmatrix} 1 \\ \lambda_3 \\ 0 \\ 0 \end{pmatrix}e^{\lambda_3\xi}+C_2\begin{pmatrix} 1 \\ \lambda_4 \\ 0 \\ 0 \end{pmatrix}e^{\lambda_4\xi}+C_3\begin{pmatrix} 0 \\ 0 \\ 1 \\ \lambda_5 \end{pmatrix}e^{\lambda_5\xi}+C_4\begin{pmatrix} 0 \\ 0 \\ 1 \\ \lambda_6 \end{pmatrix}e^{\lambda_6\xi}.$$

By applying the unstable manifold theorem, it turns out to be that, as $\xi\rightarrow-\infty$, there are $\alpha$, $\beta$, $\gamma$ and $\sigma$ such that
\begin{align*}
\phi(\xi)&=\alpha e^{\lambda_3\xi}+\beta e^{\lambda_4\xi}+h.o.t.,\\
\psi(\xi)&=\gamma e^{\lambda_5\xi}+\sigma e^{\lambda_6\xi}+h.o.t..
\end{align*}

We first note that for any $\beta<0$, no matter how large $\alpha$ is, $\phi$ is finally negative for sufficiently large negative $\xi$, since $\lambda_3>\lambda_4$. Hence, $\beta\geq 0$. Also, when $\beta=0$, $\alpha$ must be positive. Similarly, $\sigma\geqslant0$ and when $\sigma=0$, then $\gamma>0$.
\qed

Secondly, we consider only two of eigenvalues of the matrix $A_2$ are equal.
\begin{lemma}\label{lem:asym02}
When \eqref{eq:chara} has only a multiple root, then the asymptotic behavior of $(\phi,\psi)$ as $\xi\rightarrow-\infty$ is shown as follows.

$(1)$ If $\lambda_3=\lambda_4$, $\lambda_5>\lambda_6$, and $\lambda_3\neq\lambda_5$, $\lambda_6$, then
\begin{align*}
\phi(\xi)&=\alpha e^{\lambda_3\xi}-\beta\xi e^{\lambda_3\xi}+h.o.t.,\\
\psi(\xi)&=\gamma e^{\lambda_5\xi}+\sigma e^{\lambda_6\xi}+h.o.t.,
\end{align*}
where $\beta\geqslant0$, $\sigma\geqslant0$, $\alpha>0\,(\beta=0)$, $\gamma>0\,(\sigma=0)$.

$(2)$ If $\lambda_3>\lambda_4$, $\lambda_5=\lambda_6$, and $\lambda_3$, $\lambda_4\neq\lambda_5$, then
\begin{align*}
\phi(\xi)&=\alpha e^{\lambda_3\xi}+\beta e^{\lambda_4\xi}+h.o.t.,\\
\psi(\xi)&=\gamma e^{\lambda_5\xi}-\sigma\xi e^{\lambda_5\xi}+h.o.t.,
\end{align*}
where $\beta\geqslant0$, $\sigma\geqslant0$, $\alpha>0\,(\beta=0)$, $\gamma>0\,(\sigma=0)$.

$(3)$ If $\lambda_3=\lambda_5$,
$\lambda_3>\lambda_4$, $\lambda_5>\lambda_6$, and $\lambda_4\neq\lambda_6$, then
\begin{align*}
\phi(\xi)&=\alpha e^{\lambda_3\xi}+\beta e^{\lambda_4\xi}+h.o.t.,\\
\psi(\xi)&=\gamma e^{\lambda_3\xi}+\sigma e^{\lambda_6\xi}+h.o.t.,
\end{align*}
where $\beta\geqslant0$, $\sigma\geqslant0$, $\alpha>0\,(\beta=0)$, $\gamma>0\,(\sigma=0)$.

$(4)$ If $\lambda_4<\lambda_3=\lambda_6<\lambda_5$, then
\begin{align*}
\phi(\xi)&=\alpha e^{\lambda_3\xi}+\beta e^{\lambda_4\xi}+h.o.t.,\\
\psi(\xi)&=\gamma e^{\lambda_5\xi}+\sigma e^{\lambda_3\xi}+h.o.t.,
\end{align*}
where $\beta\geqslant0$, $\sigma\geqslant0$, $\alpha>0\,(\beta=0)$, $\gamma>0\,(\sigma=0)$.

$(5)$ If $\lambda_3>\lambda_4=\lambda_5>\lambda_6$, then
\begin{align*}
\phi(\xi)&=\alpha e^{\lambda_3\xi}+\beta e^{\lambda_4\xi}+h.o.t.,\\
\psi(\xi)&=\gamma e^{\lambda_4\xi}+\sigma e^{\lambda_6\xi}+h.o.t.,
\end{align*}
where $\beta\geqslant0$, $\sigma\geqslant0$, $\alpha>0\,(\beta=0)$, $\gamma>0\,(\sigma=0)$.

$(6)$ If $\lambda_4=\lambda_6$,
$\lambda_3>\lambda_4$, $\lambda_5>\lambda_6$, and $\lambda_3\neq\lambda_5$, then
\begin{align*}
\phi(\xi)&=\alpha e^{\lambda_3\xi}+\beta e^{\lambda_4\xi}+h.o.t.,\\
\psi(\xi)&=\gamma e^{\lambda_5\xi}+\sigma e^{\lambda_4\xi}+h.o.t.,
\end{align*}
where $\beta\geqslant0$, $\sigma\geqslant0$, $\alpha>0\,(\beta=0)$, $\gamma>0\,(\sigma=0)$.
\end{lemma}

\Proof First of all, consider the case (1). From $\lambda_3=\lambda_4$, we see that $c=2$. Thus $\lambda_3=\lambda_4=1$. Set
$$B_{21}:=A_2-I=\begin{pmatrix} -1 & 1 & 0 & 0 \\ -1 & 1 & 0 & 0 \\ 0 & 0 & -1 & 1 \\ 0 & 0 & -\frac{r}{d} & \frac{2}{d}-1 \end{pmatrix}.$$
From directly calculating,
$$(B_{21})^2=\begin{pmatrix} -1 & 1 & 0 & 0 \\ -1 & 1 & 0 & 0 \\ 0 & 0 & -1 & 1 \\ 0 & 0 & -\frac{r}{d} & \frac{2}{d}-1 \end{pmatrix}^2=\begin{pmatrix} 0 & 0 & 0 & 0 \\ 0 & 0 & 0 & 0 \\ 0 & 0 & 1-\frac{r}{d} & \frac{2}{d}-2 \\[0.1cm] 0 & 0 & \frac{2r}{d}-\frac{2r}{d^2} & (\frac{2}{d}-1)^2-\frac{r}{d} \end{pmatrix}.$$
Hence from the generalized characteristic equations $(B_{21})^2r=0$, we can find two linearly independent generalized eigenvectors,
$$r^2_{10}=\begin{pmatrix} 1 \\ 1 \\ 0 \\ 0 \end{pmatrix},\qquad  r^2_{20}=\begin{pmatrix} 0 \\ 1 \\ 0 \\ 0 \end{pmatrix}.$$
Thus,
$$r^2_{11}=\begin{pmatrix} -1 & 1 & 0 & 0 \\ -1 & 1 & 0 & 0 \\ 0 & 0 & -1 & 1 \\ 0 & 0 & -\frac{r}{d} & \frac{2}{d}-1 \end{pmatrix}\begin{pmatrix} 1 \\ 1 \\ 0 \\ 0 \end{pmatrix}=\begin{pmatrix} 0 \\ 0 \\ 0 \\ 0 \end{pmatrix},$$
$$r^2_{21}=\begin{pmatrix} -1 & 1 & 0 & 0 \\ -1 & 1 & 0 & 0 \\ 0 & 0 & -1 & 1 \\ 0 & 0 & -\frac{r}{d} & \frac{2}{d}-1 \end{pmatrix}\begin{pmatrix} 0 \\ 1 \\ 0 \\ 0 \end{pmatrix}=\begin{pmatrix} 1 \\ 1 \\ 0 \\ 0 \end{pmatrix}.$$
From Lemma \ref{lem:asym01}, the corresponding eigenvectors of $\lambda_5$ and $\lambda_6$ are $r^1_3$ and $r^1_4$. Consequently, there are some arbitrarily given numbers $C_1$, $C_2$, $C_3$ and $C_4$ such that
$$\begin{pmatrix} \phi_- \\ Y_- \\ \psi_- \\ Z_-  \end{pmatrix}=C_1\begin{pmatrix} 1 \\ 1 \\ 0 \\ 0 \end{pmatrix}e^{\lambda_3\xi}+C_2\left[\begin{pmatrix} 0 \\ 1 \\ 0 \\ 0 \end{pmatrix}+\begin{pmatrix} 1 \\ 1 \\ 0 \\ 0 \end{pmatrix}\xi\right]e^{\lambda_4\xi}+C_3\begin{pmatrix} 0 \\ 0 \\ 1 \\ \lambda_5 \end{pmatrix}e^{\lambda_5\xi}+C_4\begin{pmatrix} 0 \\ 0 \\ 1 \\ \lambda_6 \end{pmatrix}e^{\lambda_6\xi}.$$

From the unstable manifold theorem, it yields that, as $\xi\rightarrow-\infty$, there are $\alpha$, $\beta$, $\gamma$ and $\sigma$ such that
\begin{align*}
\phi(\xi)&=\alpha e^{\lambda_3\xi}-\beta\xi e^{\lambda_3\xi}+h.o.t.,\\
\psi(\xi)&=\gamma e^{\lambda_5\xi}+\sigma e^{\lambda_6\xi}+h.o.t..
\end{align*}
With the similar proof as in Lemma \ref{lem:asym01}, we conclude that $\beta\geqslant0$, $\sigma\geqslant0$, $\alpha>0\,(\beta=0)$, $\gamma>0\,(\sigma=0)$.

Then take the case (2) into consideration. Since $\lambda_5=\lambda_6$, then $c=2\sqrt{rd}$. Therefore, in fact, $\lambda_5=\lambda_6=\frac{c}{2d}$. Set
$$B_{22}:=A_2-\frac{c}{2d}I=\begin{pmatrix}-\frac{c}{2d} & 1 & 0 & 0 \\[0.1cm] -1 & c-\frac{c}{2d} & 0 & 0 \\[0.1cm] 0 & 0 & -\frac{c}{2d} & 1 \\[0.1cm] 0 & 0 & -\frac{r}{d} & \frac{c}{2d} \end{pmatrix}.$$
By directly calculating,
$$(B_{22})^2=\begin{pmatrix}-\frac{c}{2d} & 1 & 0 & 0 \\[0.1cm] -1 & c-\frac{c}{2d} & 0 & 0 \\[0.1cm] 0 & 0 & -\frac{c}{2d} & 1 \\[0.1cm] 0 & 0 & -\frac{r}{d} & \frac{c}{2d} \end{pmatrix}^2=\begin{pmatrix} \frac{c^2}{4d^2}-1 & c-\frac{c}{d} & 0 & 0 \\[0.1cm] \frac{c}{d}-c & (c-\frac{c}{2d})^2-1 & 0 & 0 \\[0.1cm] 0 & 0 & 0 & 0 \\[0.1cm] 0 & 0 & 0 & 0 \end{pmatrix}.$$
Thus by the generalized characteristic equations $(B_{22})^2r=0$, we can find two linearly independent generalized eigenvectors,
$$r^2_{30}=\begin{pmatrix} 0 \\ 0 \\ 1 \\ \frac{c}{2d} \end{pmatrix},\qquad r^2_{40}=\begin{pmatrix} 0 \\ 0 \\ 0 \\ 1 \end{pmatrix}.$$
Hence,
$$r^2_{31}=\begin{pmatrix}-\frac{c}{2d} & 1 & 0 & 0 \\[0.1cm] -1 & c-\frac{c}{2d} & 0 & 0 \\[0.1cm] 0 & 0 & -\frac{c}{2d} & 1 \\[0.1cm] 0 & 0 & -\frac{r}{d} & \frac{c}{2d} \end{pmatrix}\begin{pmatrix} 0 \\[0.1cm] 0 \\[0.1cm] 1 \\[0.1cm] \frac{c}{2d} \end{pmatrix}=\begin{pmatrix} 0 \\[0.1cm] 0 \\[0.1cm] 0 \\[0.1cm] 0 \end{pmatrix},$$
$$r^2_{41}=\begin{pmatrix}-\frac{c}{2d} & 1 & 0 & 0 \\[0.1cm] -1 & c-\frac{c}{2d} & 0 & 0 \\[0.1cm] 0 & 0 & -\frac{c}{2d} & 1 \\[0.1cm] 0 & 0 & -\frac{r}{d} & \frac{c}{2d} \end{pmatrix}\begin{pmatrix} 0 \\[0.1cm] 0 \\[0.1cm] 0 \\[0.1cm] 1 \end{pmatrix}=\begin{pmatrix} 0 \\[0.1cm] 0 \\[0.1cm] 1 \\[0.1cm] \frac{c}{2d} \end{pmatrix}.$$
From Lemma \ref{lem:asym01}, the corresponding eigenvectors of $\lambda_3$ and $\lambda_4$ are $r^1_1$ and $r^1_2$. As a result, there are some arbitrarily given numbers $C_1$, $C_2$, $C_3$ and $C_4$ such that
\begin{align*}
\begin{pmatrix} \phi_- \\ Y_- \\ \psi_- \\ Z_-  \end{pmatrix}=&C_1\begin{pmatrix} 1 \\ \lambda_3 \\ 0 \\ 0 \end{pmatrix}e^{\lambda_3\xi}+C_2\begin{pmatrix} 1 \\ \lambda_4 \\ 0 \\ 0 \end{pmatrix}e^{\lambda_4\xi}+C_3\begin{pmatrix} 0 \\ 0 \\ 1 \\ \frac{c}{2d} \end{pmatrix}e^{\lambda_5\xi}\\[0.2cm]
&+C_4\left[\begin{pmatrix} 0 \\ 0 \\ 0 \\ 1 \end{pmatrix}+\begin{pmatrix} 0 \\ 0 \\ 1 \\ \frac{c}{2d} \end{pmatrix}\xi\right]e^{\lambda_5\xi}.
\end{align*}

With the aid of the unstable manifold theorem and the similar proof in Lemma \ref{lem:asym01}, it is not hard to see, as $\xi\rightarrow-\infty$, there are $\alpha$, $\beta$, $\gamma$ and $\sigma$ such that
\begin{align*}
\phi(\xi)&=\alpha e^{\lambda_3\xi}+\beta e^{\lambda_4\xi}+h.o.t.,\\
\psi(\xi)&=\gamma e^{\lambda_5\xi}-\sigma\xi e^{\lambda_5\xi}+h.o.t.,
\end{align*}
where $\beta\geqslant0$, $\sigma\geqslant0$, $\alpha>0\,(\beta=0)$, $\gamma>0\,(\sigma=0)$.

Next, we will consider the case (3). In this case, due to $\lambda_3=\lambda_5$, then $d\lambda_3^2-c\lambda_3+r=0$. We also define
$$B_{23}:=A_2-\lambda_3I=\begin{pmatrix}-\lambda_3 & 1 & 0 & 0 \\ -1 & c-\lambda_3 & 0 & 0 \\0 & 0 & -\lambda_3 & 1 \\ 0 & 0 & -\frac{r}{d} & \frac{c}{d}-\lambda_3 \end{pmatrix}.$$
From directly calculating,
\begin{align*}
(B_{23})^2=&\begin{pmatrix} -\lambda_3 & 1 & 0 & 0 \\ -1 & c-\lambda_3 & 0 & 0 \\ 0 & 0 & -\lambda_3 & 1 \\ 0 & 0 & -\frac{r}{d} & \frac{c}{d}-\lambda_3 \end{pmatrix}^2\\[0.2cm]
=&\begin{pmatrix} \lambda_3^2-1 & c-2\lambda_3 & 0 & 0 \\[0.1cm] 2\lambda_3-c & (c-\lambda_3)^2-1 & 0 & 0 \\[0.1cm] 0 & 0 & \lambda_3^2-\frac{r}{d} & \frac{c}{d}-2\lambda_3 \\[0.1cm] 0 & 0 & \frac{2r}{d}\lambda_3-\frac{rc}{d^2} & (\frac{c}{d}-\lambda_3)^2-\frac{r}{d} \end{pmatrix}.
\end{align*}
Therefore, by the generalized characteristic equations $(B_{23})^2r=0$, we can also find two linearly independent generalized eigenvectors,
$$r^2_{50}=\begin{pmatrix} 1 \\ \lambda_3 \\ 0 \\ 0 \end{pmatrix},\qquad r^2_{60}=\begin{pmatrix} 0 \\ 0 \\ 1 \\ \lambda_3 \end{pmatrix}.$$
Thus,
$$r^2_{51}=\begin{pmatrix}-\lambda_3 & 1 & 0 & 0 \\ -1 & c-\lambda_3 & 0 & 0 \\0 & 0 & -\lambda_3 & 1 \\ 0 & 0 & -\frac{r}{d} & \frac{c}{d}-\lambda_3 \end{pmatrix}\begin{pmatrix} 1 \\ \lambda_3 \\ 0 \\ 0 \end{pmatrix}=\begin{pmatrix} 0 \\ 0 \\ 0 \\ 0 \end{pmatrix},$$
$$r^2_{61}=\begin{pmatrix}-\lambda_3 & 1 & 0 & 0 \\ -1 & c-\lambda_3 & 0 & 0 \\0 & 0 & -\lambda_3 & 1 \\ 0 & 0 & -\frac{r}{d} & \frac{c}{d}-\lambda_3 \end{pmatrix}\begin{pmatrix} 0 \\ 0 \\ 1 \\ \lambda_3 \end{pmatrix}=\begin{pmatrix} 0 \\ 0 \\ 0 \\ 0 \end{pmatrix}.$$
Similarly note that from Lemma \ref{lem:asym01}, the corresponding eigenvectors of $\lambda_4$ and $\lambda_6$ are $r^1_2$ and $r^1_4$. Hence, there are some arbitrarily given numbers $C_1$, $C_2$, $C_3$ and $C_4$ such that
\begin{align*}
\begin{pmatrix} \phi_- \\ Y_- \\ \psi_- \\ Z_-  \end{pmatrix}=&C_1\begin{pmatrix} 1 \\ \lambda_3 \\ 0 \\ 0 \end{pmatrix}e^{\lambda_3\xi}+C_2\begin{pmatrix} 1 \\ \lambda_4 \\ 0 \\ 0 \end{pmatrix}e^{\lambda_4\xi}+C_3\begin{pmatrix} 0 \\ 0 \\ 1 \\ \lambda_3 \end{pmatrix}e^{\lambda_3\xi}\\[0.2cm]
&+C_4\begin{pmatrix} 0 \\ 0 \\ 1 \\ \lambda_6 \end{pmatrix}e^{\lambda_6\xi}.
\end{align*}

Similar, we can conclude that
\begin{align*}
\phi(\xi)&=\alpha e^{\lambda_3\xi}+\beta e^{\lambda_4\xi}+h.o.t.,\\
\psi(\xi)&=\gamma e^{\lambda_3\xi}+\sigma e^{\lambda_6\xi}+h.o.t.,
\end{align*}
where $\beta\geqslant0$, $\sigma\geqslant0$, $\alpha>0\,(\beta=0)$, $\gamma>0\,(\sigma=0)$.

The rest of the proof is similar to the proof of the case (3).
\qed

Thirdly, the matrix $A_2$ having two different multiple eigenvalues is taken into consideration, in fact, that is $\lambda_3=\lambda_4$, $\lambda_5=\lambda_6$, $\lambda_3\neq\lambda_5$, or $\lambda_3=\lambda_5$, $\lambda_4=\lambda_6$, $\lambda_3>\lambda_4$.

\begin{lemma}\label{lem:asym03}
When \eqref{eq:chara} has two different multiple roots, then the asymptotic behavior of $(\phi,\psi)$ as $\xi\rightarrow-\infty$ is shown in the following.

$(1)$ If $\lambda_3=\lambda_4$, $\lambda_5=\lambda_6$, and $\lambda_3\neq\lambda_5$, then
\begin{align*}
\phi(\xi)&=\alpha e^{\lambda_3\xi}-\beta\xi e^{\lambda_3\xi}+h.o.t.,\\
\psi(\xi)&=\gamma e^{\lambda_5\xi}-\sigma\xi e^{\lambda_5\xi}+h.o.t.,
\end{align*}
where $\beta\geqslant0$, $\sigma\geqslant0$, $\alpha>0\,(\beta=0)$, $\gamma>0\,(\sigma=0)$.

$(2)$ If $\lambda_3=\lambda_5$, $\lambda_4=\lambda_6$, and $\lambda_3>\lambda_4$, then
\begin{align*}
\phi(\xi)&=\alpha e^{\lambda_3\xi}+\beta e^{\lambda_4\xi}+h.o.t.,\\
\psi(\xi)&=\gamma e^{\lambda_3\xi}+\sigma e^{\lambda_4\xi}+h.o.t.,
\end{align*}
where $\beta\geqslant0$, $\sigma\geqslant0$, $\alpha>0\,(\beta=0)$, $\gamma>0\,(\sigma=0)$.
\end{lemma}

\Proof From the proof of the case (1) and (2) in Lemma \ref{lem:asym02}, it turns out that there are some arbitrarily given numbers $C_1$, $C_2$, $C_3$ and $C_4$ such that
\begin{align*}
\begin{pmatrix} \phi_- \\ Y_- \\ \psi_- \\ Z_-  \end{pmatrix}=&C_1\begin{pmatrix} 1 \\ 1 \\ 0 \\ 0 \end{pmatrix}e^{\lambda_3\xi}+C_2\left[\begin{pmatrix} 0 \\ 1 \\ 0 \\ 0 \end{pmatrix}+\begin{pmatrix} 1 \\ 1 \\ 0 \\ 0 \end{pmatrix}\xi\right]e^{\lambda_3\xi}+C_3\begin{pmatrix} 0 \\ 0 \\ 1 \\ \frac{c}{2d} \end{pmatrix}e^{\lambda_5\xi}\\[0.2cm]
&+C_4\left[\begin{pmatrix} 0 \\ 0 \\ 0 \\ 1 \end{pmatrix}+\begin{pmatrix} 0 \\ 0 \\ 1 \\ \frac{c}{2d} \end{pmatrix}\xi\right]e^{\lambda_5\xi},
\end{align*}
which leads to the desired result followed from the unstable manifold theorem.

Similarly, from the proof of the case (3) and (6) in Lemma \ref{lem:asym02} and with the aid of the unstable manifold theorem, the rest of the result can be obtained.
\qed

Fourthly, we will research the asymptotic behavior, when the matrix $A_2$ has a triple eigenvalue.

\begin{lemma}\label{lem:asym04}
When \eqref{eq:chara} has a triple root, then the asymptotic behavior of $(\phi,\psi)$ as $\xi\rightarrow-\infty$ is shown as follows.

$(1)$ If $\lambda_3=\lambda_4=\lambda_5>\lambda_6$, then
\begin{align*}
\phi(\xi)&=\alpha e^{\lambda_3\xi}-\beta\xi e^{\lambda_3\xi}+h.o.t.,\\
\psi(\xi)&=\gamma e^{\lambda_3\xi}+\sigma e^{\lambda_6\xi}+h.o.t.,
\end{align*}
where $\beta\geqslant0$, $\sigma\geqslant0$, $\alpha>0\,(\beta=0)$, $\gamma>0\,(\sigma=0)$.

$(2)$ If $\lambda_5>\lambda_6=\lambda_3=\lambda_4$, then
\begin{align*}
\phi(\xi)&=\alpha e^{\lambda_3\xi}-\beta\xi e^{\lambda_3\xi}+h.o.t.,\\
\psi(\xi)&=\gamma e^{\lambda_5\xi}+\sigma e^{\lambda_3\xi}+h.o.t.,
\end{align*}
where $\beta\geqslant0$, $\sigma\geqslant0$, $\alpha>0\,(\beta=0)$, $\gamma>0\,(\sigma=0)$.

$(3)$ If $\lambda_5=\lambda_6=\lambda_3>\lambda_4$, then
\begin{align*}
\phi(\xi)&=\alpha e^{\lambda_3\xi}+\beta e^{\lambda_4\xi}+h.o.t.,\\
\psi(\xi)&=\gamma e^{\lambda_3\xi}-\sigma\xi e^{\lambda_3\xi}+h.o.t.,
\end{align*}
where $\beta\geqslant0$, $\sigma\geqslant0$, $\alpha>0\,(\beta=0)$, $\gamma>0\,(\sigma=0)$.

$(4)$ If $\lambda_3>\lambda_4=\lambda_5=\lambda_6$, then
\begin{align*}
\phi(\xi)&=\alpha e^{\lambda_3\xi}+\beta e^{\lambda_4\xi}+h.o.t.,\\
\psi(\xi)&=\gamma e^{\lambda_4\xi}-\sigma\xi e^{\lambda_4\xi}+h.o.t.,
\end{align*}
where $\beta\geqslant0$, $\sigma\geqslant0$, $\alpha>0\,(\beta=0)$, $\gamma>0\,(\sigma=0)$.
\end{lemma}

\Proof We firstly consider the case (1), namely, $\lambda_3=\lambda_4=\lambda_5>\lambda_6$. From the proof in Lemma \ref{lem:asym02}, $1=\lambda_3=\lambda_4=\lambda_5$. Thus $1-\frac{2}{d}+\frac{r}{d}=0$, which leads to $d+r=2$. Then it is not hard to calculate
\begin{align*}
(B_{21})^3=&\begin{pmatrix} -1 & 1 & 0 & 0 \\ -1 & 1 & 0 & 0 \\ 0 & 0 & -1 & 1 \\ 0 & 0 & -\frac{r}{d} & \frac{2}{d}-1 \end{pmatrix}^3\\[0.2cm]
=&\begin{pmatrix} 0 & 0 & 0 & 0 \\ 0 & 0 & 0 & 0 \\[0.1cm] 0 & 0 & -1+\frac{3r}{d}-\frac{2r}{d^2} & 3-\frac{6}{d}-\frac{r}{d}+\frac{4}{d^2} \\[0.1cm] 0 & 0 & -\frac{r}{d}+\frac{r^2}{d^2}+(\frac{2}{d}-1)(\frac{2r}{d}-\frac{2r}{d^2}) & -\frac{2r}{d^2}+\frac{2r}{d}+(\frac{2}{d}-1)^3-\frac{r}{d}(\frac{2}{d}-1) \end{pmatrix}.
\end{align*}
Hence, by the generalized characteristic equations $(B_{21})^3r=0$, we can find three linearly independent generalized eigenvectors,
$$r^3_{10}=\begin{pmatrix} 1 \\ 1 \\ 0 \\ 0 \end{pmatrix},\qquad r^3_{20}=\begin{pmatrix} 0 \\ 1 \\ 0 \\ 0 \end{pmatrix},\qquad  r^3_{30}=\begin{pmatrix} 0 \\ 0 \\ 1 \\ 1 \end{pmatrix}.$$
Thus,
\begin{align*}
r^3_{11}=\begin{pmatrix} -1 & 1 & 0 & 0 \\ -1 & 1 & 0 & 0 \\0 & 0 & -1 & 1 \\ 0 & 0 & -\frac{r}{d} & \frac{2}{d}-1 \end{pmatrix}\begin{pmatrix} 1 \\ 1 \\ 0 \\ 0 \end{pmatrix}=\begin{pmatrix} 0 \\ 0 \\ 0 \\ 0 \end{pmatrix},\\[0.2cm]
r^3_{21}=\begin{pmatrix} -1 & 1 & 0 & 0 \\ -1 & 1 & 0 & 0 \\0 & 0 & -1 & 1 \\ 0 & 0 & -\frac{r}{d} & \frac{2}{d}-1 \end{pmatrix}\begin{pmatrix} 0 \\ 1 \\ 0 \\ 0 \end{pmatrix}=\begin{pmatrix} 1 \\ 1 \\ 0 \\ 0 \end{pmatrix},\\[0.2cm]
r^3_{22}=\begin{pmatrix} -1 & 1 & 0 & 0 \\ -1 & 1 & 0 & 0 \\0 & 0 & -1 & 1 \\ 0 & 0 & -\frac{r}{d} & \frac{2}{d}-1 \end{pmatrix}\begin{pmatrix} 1 \\ 1 \\ 0 \\ 0 \end{pmatrix}=\begin{pmatrix} 0 \\ 0 \\ 0 \\ 0 \end{pmatrix},\\[0.2cm]
r^3_{31}=\begin{pmatrix} -1 & 1 & 0 & 0 \\ -1 & 1 & 0 & 0 \\0 & 0 & -1 & 1 \\ 0 & 0 & -\frac{r}{d} & \frac{2}{d}-1 \end{pmatrix}\begin{pmatrix} 0 \\ 0 \\ 1 \\ 1 \end{pmatrix}=\begin{pmatrix} 0 \\ 0 \\ 0 \\ 0 \end{pmatrix}.\\[0.2cm]
\end{align*}
From Lemma \ref{lem:asym02}, the corresponding eigenvector of $\lambda_6$ is $r^1_4$. Therefore, there are some arbitrarily given numbers $C_1$, $C_2$, $C_3$ and $C_4$ such that
$$\begin{pmatrix} \phi_- \\ Y_- \\ \psi_- \\ Z_-  \end{pmatrix}=C_1\begin{pmatrix} 1 \\ 1 \\ 0 \\ 0 \end{pmatrix}e^{\lambda_3\xi}+C_2\left[\begin{pmatrix} 0 \\ 1 \\ 0 \\ 0 \end{pmatrix}+\begin{pmatrix} 1 \\ 1 \\ 0 \\ 0 \end{pmatrix}\xi\right]e^{\lambda_3\xi}+C_3\begin{pmatrix} 0 \\ 0 \\ 1 \\ 1 \end{pmatrix}e^{\lambda_3\xi}+C_4\begin{pmatrix} 0 \\ 0 \\ 1 \\ \lambda_6 \end{pmatrix}e^{\lambda_6\xi}.$$
Similarly, we can obtain
\begin{align*}
\phi(\xi)&=\alpha e^{\lambda_3\xi}-\beta\xi e^{\lambda_3\xi}+h.o.t.,\\
\psi(\xi)&=\gamma e^{\lambda_3\xi}+\sigma e^{\lambda_6\xi}+h.o.t.,
\end{align*}
where $\beta\geqslant0$, $\sigma\geqslant0$, $\alpha>0\,(\beta=0)$, $\gamma>0\,(\sigma=0)$.

We then consider the case (3). From the proof in Lemma \ref{lem:asym02}, easily calculating
\begin{align*}
(B_{22})^3=&\begin{pmatrix}-\frac{c}{2d} & 1 & 0 & 0 \\ -1 & c-\frac{c}{2d} & 0 & 0 \\0 & 0 & -\frac{c}{2d} & 1 \\ 0 & 0 & -\frac{r}{d} & \frac{c}{2d} \end{pmatrix}^3\\[0.2cm]
=&\begin{pmatrix} -\frac{c^3}{8d^3}+\frac{3c}{2d}-c & \frac{3c^2}{4d^2}-\frac{3c^2}{2d}+c^2-1 & 0 & 0 \\[0.1cm] -\frac{3c^2}{4d^2}+\frac{3c^2}{2d}-c^2+1 & \frac{c}{d}-c+(c-\frac{c}{2d})^3-(c-\frac{c}{2d}) & 0 & 0 \\[0.1cm]  0 & 0 & 0 & 0 \\ 0 & 0 & 0 & 0 \end{pmatrix}.
\end{align*}
Since $\frac{c}{2d}=\lambda_5=\lambda_6=\lambda_3$, then $(\frac{c}{2d})^2-\frac{c^2}{2d}+1=0$, that is $\frac{c^2}{4d^2}=\frac{c^2}{2d}-1$. With  $\frac{c^2}{4d^2}=\frac{c^2}{2d}-1$, we can simplify $(B_{22})^3$, which is
$$\begin{pmatrix} -\frac{c^3}{2d}+\frac{2c}{d} & c^2-4 & 0 & 0 \\[0.2cm] -c^2+4 & c^3-\frac{c^3}{2d}-4c+\frac{2c}{d} & 0 & 0 \\ 0 & 0 & 0 & 0 \\ 0 & 0 & 0 & 0 \end{pmatrix}.$$
Therefore, by the generalized characteristic equations $(B_{22})^3r=0$, we can also find three linearly independent generalized eigenvectors,
$$r^3_{40}=\begin{pmatrix} 1 \\ \frac{c}{2d} \\ 0 \\ 0 \end{pmatrix},\qquad
r^3_{50}=\begin{pmatrix} 0 \\ 0 \\ 1 \\ \frac{c}{2d} \end{pmatrix},\qquad r^3_{60}=\begin{pmatrix} 0 \\ 0 \\ 0 \\ 1 \end{pmatrix}.$$
Thus,
\begin{align*}
r^3_{41}=\begin{pmatrix}-\frac{c}{2d} & 1 & 0 & 0 \\ -1 & c-\frac{c}{2d} & 0 & 0 \\0 & 0 & -\frac{c}{2d} & 1 \\ 0 & 0 & -\frac{r}{d} & \frac{c}{2d} \end{pmatrix}\begin{pmatrix} 1 \\ \frac{c}{2d} \\ 0 \\ 0 \end{pmatrix}=\begin{pmatrix} 0 \\ 0 \\ 0 \\ 0 \end{pmatrix},\\[0.2cm]
r^3_{51}=\begin{pmatrix}-\frac{c}{2d} & 1 & 0 & 0 \\ -1 & c-\frac{c}{2d} & 0 & 0 \\0 & 0 & -\frac{c}{2d} & 1 \\ 0 & 0 & -\frac{r}{d} & \frac{c}{2d} \end{pmatrix}\begin{pmatrix} 0 \\ 0 \\ 1 \\ \frac{c}{2d} \end{pmatrix}=\begin{pmatrix} 0 \\ 0 \\ 0 \\ 0 \end{pmatrix},\\[0.2cm]
r^3_{61}=\begin{pmatrix}-\frac{c}{2d} & 1 & 0 & 0 \\ -1 & c-\frac{c}{2d} & 0 & 0 \\0 & 0 & -\frac{c}{2d} & 1 \\ 0 & 0 & -\frac{r}{d} & \frac{c}{2d} \end{pmatrix}\begin{pmatrix} 0 \\ 0 \\ 0 \\ 1 \end{pmatrix}=\begin{pmatrix} 0 \\ 0 \\ 1 \\ \frac{c}{2d} \end{pmatrix},\\[0.2cm]
r^3_{62}=\begin{pmatrix} -\frac{c}{2d} & 1 & 0 & 0 \\ -1 & c-\frac{c}{2d} & 0 & 0 \\ 0 & 0 & -\frac{c}{2d} & 1 \\ 0 & 0 & -\frac{r}{d} & \frac{c}{2d} \end{pmatrix}\begin{pmatrix} 0 \\ 0 \\ 1 \\ \frac{c}{2d} \end{pmatrix}=\begin{pmatrix} 0 \\ 0 \\ 0 \\ 0 \end{pmatrix}.\\[0.2cm]
\end{align*}
Also from Lemma \ref{lem:asym02}, the corresponding eigenvector of $\lambda_4$ is $r^1_2$. Therefore, there are some arbitrarily given numbers $C_1$, $C_2$, $C_3$ and $C_4$ such that
\begin{align*}
\begin{pmatrix} \phi_- \\ Y_- \\ \psi_- \\ Z_-  \end{pmatrix}=&C_1\begin{pmatrix} 1 \\ \frac{c}{2d} \\ 0 \\ 0 \end{pmatrix}e^{\lambda_3\xi}+C_2\begin{pmatrix} 1 \\ \lambda_4 \\ 0 \\ 0 \end{pmatrix}e^{\lambda_4\xi}+C_3\begin{pmatrix} 0 \\ 0 \\ 1 \\ \frac{c}{2d} \end{pmatrix}e^{\lambda_3\xi}\\[0.2cm]
&+C_4\left[\begin{pmatrix} 0 \\ 0 \\ 0 \\ 1 \end{pmatrix}+\begin{pmatrix} 0 \\ 0 \\ 1 \\ \frac{c}{2d} \end{pmatrix}\xi\right]e^{\lambda_3\xi}.
\end{align*}
The rest of the proof is similar.
\qed

In the end, we will consider all the eigenvalues of the matrix $A_2$ are equal, that is $\lambda_3=\lambda_4=\lambda_5=\lambda_6$.

\begin{lemma}\label{lem:asym05}
If $\lambda_3=\lambda_4=\lambda_5=\lambda_6$, then the asymptotic behavior of $(\phi,\psi)$ as $\xi\rightarrow-\infty$ is shown in the following
\begin{align*}
\phi(\xi)&=\alpha e^{\lambda_3\xi}-\beta\xi e^{\lambda_3\xi}+h.o.t.,\\
\psi(\xi)&=\gamma e^{\lambda_3\xi}-\sigma\xi e^{\lambda_3\xi}+h.o.t.,
\end{align*}
where $\beta\geqslant0$, $\sigma\geqslant0$, $\alpha>0\,(\beta=0)$, $\gamma>0\,(\sigma=0)$.
\end{lemma}

\Proof Since $\lambda_3=\lambda_4=\lambda_5=\lambda_6$, then all them are equal to 1 and also $d=r=1$. Set
$$B_{41}:=A_2-I=\begin{pmatrix} -1 & 1 & 0 & 0 \\ -1 & 1 & 0 & 0 \\ 0 & 0 & -1 & 1 \\ 0 & 0 & -1 & 1 \end{pmatrix}.$$
and it is not hard to calculate,
$$(B_{41})^2=(B_{41})^3=(B_{41})^4=\begin{pmatrix} 0 & 0 & 0 & 0 \\ 0 & 0 & 0 & 0 \\ 0 & 0 & 0 & 0 \\ 0 & 0 & 0 & 0 \end{pmatrix}.$$
Therefore from the generalized characteristic equations $(B_{41})^4r=0$, we can find four linearly independent generalized eigenvectors,
$$r^4_{10}=\begin{pmatrix} 1 \\ 1 \\ 0 \\ 0 \end{pmatrix},\qquad  r^4_{20}=\begin{pmatrix} 0 \\ 1 \\ 0 \\ 0 \end{pmatrix},\qquad
r^4_{30}=\begin{pmatrix} 0 \\ 0 \\ 1 \\ 1 \end{pmatrix},\qquad  r^4_{40}=\begin{pmatrix} 0 \\ 0 \\ 0 \\ 1 \end{pmatrix}.$$
Thus,
\begin{align*}
r^4_{11}=\begin{pmatrix} -1 & 1 & 0 & 0 \\ -1 & 1 & 0 & 0 \\0 & 0 & -1 & 1 \\ 0 & 0 & -1 & 1 \end{pmatrix}\begin{pmatrix} 1 \\ 1 \\ 0 \\ 0 \end{pmatrix}=\begin{pmatrix} 0 \\ 0 \\ 0 \\ 0 \end{pmatrix},\\[0.2cm]
r^4_{21}=\begin{pmatrix} -1 & 1 & 0 & 0 \\ -1 & 1 & 0 & 0 \\0 & 0 & -1 & 1 \\ 0 & 0 & -1 & 1 \end{pmatrix}\begin{pmatrix} 0 \\ 1 \\ 0 \\ 0 \end{pmatrix}=\begin{pmatrix} 1 \\ 1 \\ 0 \\ 0 \end{pmatrix},\\[0.2cm]
r^4_{22}=\begin{pmatrix} -1 & 1 & 0 & 0 \\ -1 & 1 & 0 & 0 \\0 & 0 & -1 & 1 \\ 0 & 0 & -1 & 1 \end{pmatrix}\begin{pmatrix} 1 \\ 1 \\ 0 \\ 0 \end{pmatrix}=\begin{pmatrix} 0 \\ 0 \\ 0 \\ 0 \end{pmatrix},\\[0.2cm]
r^4_{31}=\begin{pmatrix} -1 & 1 & 0 & 0 \\ -1 & 1 & 0 & 0 \\0 & 0 & -1 & 1 \\ 0 & 0 & -1 & 1 \end{pmatrix}\begin{pmatrix} 0 \\ 0 \\ 1 \\ 1 \end{pmatrix}=\begin{pmatrix} 0 \\ 0 \\ 0 \\ 0 \end{pmatrix},\\[0.2cm]
r^4_{41}=\begin{pmatrix} -1 & 1 & 0 & 0 \\ -1 & 1 & 0 & 0 \\0 & 0 & -1 & 1 \\ 0 & 0 & -1 & 1 \end{pmatrix}\begin{pmatrix} 0 \\ 0 \\ 0 \\ 1 \end{pmatrix}=\begin{pmatrix} 0 \\ 0 \\ 1 \\ 1 \end{pmatrix},\\[0.2cm]
r^4_{42}=\begin{pmatrix} -1 & 1 & 0 & 0 \\ -1 & 1 & 0 & 0 \\0 & 0 & -1 & 1 \\ 0 & 0 & -1 & 1 \end{pmatrix}\begin{pmatrix} 0 \\ 0 \\ 1 \\ 1 \end{pmatrix}=\begin{pmatrix} 0 \\ 0 \\ 0 \\ 0 \end{pmatrix}.
\end{align*}
Consequently, there are some arbitrarily given numbers $C_1$, $C_2$, $C_3$ and $C_4$ such that
\begin{align*}
\begin{pmatrix}\phi_- \\ Y_- \\ \psi_- \\ Z_-  \end{pmatrix}=&C_1\begin{pmatrix} 1 \\ 1 \\ 0 \\ 0 \end{pmatrix}e^{\lambda_3\xi}+C_2\left[\begin{pmatrix} 0 \\ 1 \\ 0 \\ 0 \end{pmatrix}+\begin{pmatrix} 1 \\ 1 \\ 0 \\ 0 \end{pmatrix}\xi\right]e^{\lambda_3\xi}+C_3\begin{pmatrix} 0 \\ 0 \\ 1 \\ 1 \end{pmatrix}e^{\lambda_3\xi}\\[0.2cm]
&+C_4\left[\begin{pmatrix} 0 \\ 0 \\ 0 \\ 1 \end{pmatrix}+\begin{pmatrix} 0 \\ 0 \\ 1 \\ 1 \end{pmatrix}\xi\right]e^{\lambda_3\xi}.
\end{align*}

From the unstable manifold theorem, it yields that, as $\xi\rightarrow-\infty$, there are $\alpha$, $\beta$, $\gamma$ and $\sigma$ such that
\begin{align*}
\phi(\xi)&=\alpha e^{\lambda_3\xi}-\beta\xi e^{\lambda_3\xi}+h.o.t.,\\
\psi(\xi)&=\gamma e^{\lambda_3\xi}-\sigma\xi e^{\lambda_3\xi}+h.o.t..
\end{align*}
With the similar proof in Lemma \ref{lem:asym01}, we conclude that $\beta\geqslant0$, $\sigma\geqslant0$, $\alpha>0\,(\beta=0)$, $\gamma>0\,(\sigma=0)$.
\qed

We end this section with three remarks, two of which are some sufficient and necessary conditions for $\lambda_3=\lambda_4=\lambda_5$, $\lambda_3=\lambda_4=\lambda_6$, $\lambda_3=\lambda_5=\lambda_6$, $\lambda_4=\lambda_5=\lambda_6$ and $\lambda_3=\lambda_4=\lambda_5=\lambda_6$, while the rest one is devoted to comparing the asymptotic behavior discussed as above with some known results.

\begin{remark}\label{remark1}
The sufficient and necessary condition under which $\lambda_3=\lambda_4=\lambda_5>\lambda_6$ holds is $c=2$ and $d=2-r>1$, and $\lambda_3=\lambda_4=\lambda_6<\lambda_5$ if and only if $c=2$ and $0<d=2-r<1$. The remaining case is $c=2$ and $d=r=1$, which leads to $\lambda_3=\lambda_4=\lambda_5=\lambda_6$. Thus, if $c=2$,$d\neq2-r$, then $\lambda_3=\lambda_4\not =\lambda_5$, $\lambda_3=\lambda_4\not =\lambda_6$.
\end{remark}

\begin{remark}\label{remark2}
First, $\lambda_5=\lambda_6$ if and only if $c=2\sqrt{rd}$, which implies that $rd\geq 1$.Thus, $\lambda_5=\lambda_6=\lambda_3>\lambda_4$ if and only if $c=2\sqrt{rd}$ and $0<d=\frac{r}{2r-1}<1$, $\lambda_5=\lambda_6=\lambda_4<\lambda_3$ if and only if $c=2\sqrt{rd}$ and $d=\frac{r}{2r-1}>1$,  $\lambda_5=\lambda_6=\lambda_3=\lambda_4$ if and only if $c=2$ and $d=r=1$. Therefore, if $c=2\sqrt{rd}$ and $d\neq\frac{r}{2r-1}$ then  $\lambda_5=\lambda_6\neq\lambda_3$, $\lambda_5=\lambda_6\neq\lambda_4$.
\end{remark}

\begin{remark}\label{remark3}
In the paper \cite{l08}, the author had already obtained the asymptotic behavior of the solution of \eqref{eq:OLV}-\eqref{eq:MC} as $\xi\rightarrow-\infty$ under the case $(iv)$ and $c>c_{\min}$, that is,
$$\lim\limits_{\xi\rightarrow-\infty}\Bigl(\phi(\xi)e^{-\tilde{\lambda}_1(\xi+\tilde{h}_1)},\,
\psi(\xi)e^{-\tilde{\lambda}_2(\xi+\tilde{h}_2)}\Bigr)=(1,1)$$
for some $\tilde{\lambda}_i$, $\tilde{h}_i>0$, $i=1,2$. Here, under the case $(iv)$ and $c\geqslant c_{\min}$, our results contain the asymptotic behavior of the solution of \eqref{eq:OLV}-\eqref{eq:MC} as $\xi\rightarrow\pm\infty$. We also remark that in \cite{wl13}, the authors stated the asymptotic behavior of traveling front solutions for the classical Lotka-Volterra cooperative system as $\xi\rightarrow\pm\infty$, which connects the origin and the positive equilibrium. Here, we carefully analyze the asymptotic behavior of traveling front solutions for Lotka-Volterra competitive system as $\xi\rightarrow-\infty$, and accurately show that only the smaller negative eigenvalue of linearization matrix of \eqref{eq:OLVS} at $(u^{\ast},0,v^{\ast},0)$ plays role in the asymptotic behavior as $\xi\rightarrow+\infty$.
\end{remark}

\section{Existence of entire solutions via a pair of coupled super-sub solutions}
In this section, we introduce the definition of a pair of coupled super-sub solutions of general quasi-monotone decreasing reaction diffusion systems, and also state the existence of entire solutions, which has been already given in \cite{cg05}, \cite{mt09}, \cite{wl15}. Here we will investigate the qualitative properties of these entire solutions such as the monotonicity, the symmetry.

Now, we focus on the following general reaction diffusion system
\begin{equation}\label{eq:gen}
\left\{\begin{array}{ll}
\partial_tu(x,t)=\partial_{xx}u(x,t)+f(u,v),\\[0.2cm]
\partial_tv(x,t)=d\,\partial_{xx}v(x,t)+g(u,v),
\end{array}
\right.
\end{equation}
where $d>0$, the reaction terms $f$ and $g$ are quasi-monotone decreasing in the following sense $\partial_{v}f(u,v)\leqslant0$, $\partial_{u}g(u,v)\leqslant0$.

First of all, we estimate the derivatives of the solutions to \eqref{eq:gen} in the following theorem, which should be in the proof of the existence and the qualitative properties of entire solutions. Although the theorem has been proved in \cite{wl15}, we restate it here for the reader's convenience.
\begin{theorem}\label{thm:bs}
For $t>1$, assume that $(u,v)$ satisfies $(0,0)\leqslant(u,v)\leqslant(K,K)$ for some positive constant $K$ and is a solution of system \eqref{eq:gen} with continuous and bounded initial functions. We also suppose that the functions $f$ and $g$ are continuous on $[0,K]\times[0,K]$ as well. Then, for $t>1$, $|\partial_tu|$, $|\partial_tv|$, $|\partial_xu|$, $|\partial_xv|$, $|\partial_{xx}u|$, $|\partial_{xx}v|$ are bounded.
\end{theorem}

Secondly, we introduce the definition of a pair of coupled super-sub solutions of \eqref{eq:gen}.
\begin{definition}\label{def:def1}
For $(x,t)\in\mathbb{R}\times(T_1,T_0)$, where $T_1<T_0\in\mathbb{R}\cup\{+\infty\}$, let $\overline u(x,t)$, $\overline v(x,t)$, $\underline u(x,t)$ and $\underline v(x,t)$ be smooth functions satisfying $(\overline u,\overline v)\geqslant(\underline u,\underline v)$ and
\begin{equation}\label{condi:csl}
\left\{\begin{array}{ll}
\partial_t{\overline u}-\partial_{xx}{\overline u}
-f(\overline u,\underline v)\geqslant0
\geqslant\partial_t{\underline u}-\partial_{xx}{\underline u}
-f(\underline u,\overline v),
\\[0.2cm]
\partial_t{\overline v}-d\,\partial_{xx}{\overline v}
-g(\underline u,\overline v)\geqslant0
\geqslant\partial_t{\underline v}-d\,\partial_{xx}{\underline v}
-g(\overline u,\underline v),
\end{array}
\right.
\end{equation}
then the pair of the functions $(\overline u(x,t),\overline v(x,t))$ and $(\underline u(x,t),\underline v(x,t))$ are called a pair of coupled super-sub solutions of \eqref{eq:gen} on $\mathbb{R}\times(T_1,T_0)$. If for any $T_1<\min\{T_0,0\}$, if $(\overline u(x,t),\overline v(x,t))$ and $(\underline u(x,t),\underline v(x,t))$ are the pair of coupled super-sub solutions of \eqref{eq:gen} on $\mathbb{R}\times(T_1,T_0)$, then we call $(\overline u(x,t),\overline v(x,t))$ and $(\underline u(x,t),\underline v(x,t))$ the pair of coupled super-sub solutions of \eqref{eq:gen} on $\mathbb{R}\times(-\infty,T_0)$. The pair of coupled super-sub solutions is called deterministic via translation, if there exist functions $\rho^u_i(t)$ and $\rho^v_i(t)$ $(i=1,2)$, such that for $(x,t)\in\mathbb{R}\times(-\infty,T_0)$
\begin{align*}
&(\overline u(x,t),\overline v(x,t))\leqslant(\underline u(x+\rho^u_1(t),t+\rho^u_2(t)),\underline v(x+\rho^v_1(t),t+\rho^v_2(t))),\\[0.2cm]
&\lim\limits_{t\rightarrow-\infty}\{|\rho^u_1(t)|+|\rho^v_1(t)|
+|\rho^u_2(t)|+|\rho^v_2(t)|\}=0.
\end{align*}
\end{definition}

Finally, the existence of entire solutions of \eqref{eq:gen} can be established from some suitable pairs of coupled super-sub solutions defined in Definition \ref{def:def1}, and the estimates of the derivatives of the solutions to \eqref{eq:gen}.
\begin{theorem}\label{thm:exig}
If $(\overline{u}(x,t),\overline{v}(x,t))$ and $(\underline{u}(x,t),\underline{v}(x,t))$ are a pair of coupled super-sub solutions of \eqref{eq:gen} on $\mathbb{R}\times(-\infty,T_0)$ for some $T_0\in\mathbb{R}\cup\{+\infty\}$, satisfying $(0,0)\leqslant(\underline{u}(x,t),\underline{v}(x,t))\leqslant(\overline{u}(x,t),
\overline{v}(x,t))\leqslant(K,K)$ for all $(x,t)\in \mathbb{R}\times(-\infty,T_0)$, where $K$ is a positive constant. Then \eqref{eq:gen} admits an entire solution $(u(x,t),v(x,t))$ such that $(\underline{u}(x,t),\underline{v}(x,t))
\leqslant(u(x,t),v(x,t))\leqslant(\overline{u}(x,t),\overline{v}(x,t))$ for all $(x,t)\in\mathbb{R}\times(-\infty,T_0)$. This entire solution is unique if the pair of coupled super-sub solutions is deterministic via translation. In addition, if $\underline u(x,t)=\underline u(-x,t)$ and $\underline v(x,t)=\underline v(-x,t)$, then $u(x,t)=u(-x,t)$ and $v(x,t)=v(-x,t)$.
\end{theorem}

\Proof By Theorem 2 in \cite{t82} and refining the proof in \cite{cg05}, we can get the existence and uniqueness of entire solutions. Now we suppose that $(u_n(x,t),v_n(x,t))$ is the unique classical bounded solution of \eqref{eq:gen} with $(u_n(x,-n),v_n(x,-n))=(\underline u(x,-n),\underline v(x,-n))$ for $x\in\mathbb{R}$ and $t>-n$. By the boundedness of $\partial_tu_n$, $\partial_tv_n$, $\partial_xu_n$, $\partial_xv_n$, $\partial_{xx}u_n$, $\partial_{xx}v_n$, which are gotten from Theorem \ref{thm:bs}, there are subsequences of $u_n(x,t)$ and $v_n(x,t)$ which converge to $u(x,t)$ and $v(x,t)$ in $C^{2,1}_{loc}(\mathbb{R}\times(-\infty,T_0))$, respectively. By the above convergence and remarking $(\underline u(x,-n),\underline v(x,-n))=(\underline u(-x,-n),\underline v(-x,-n))$, we can get the symmetry.
\qed
\begin{remark}\label{remark4}
If $\overline u(x,t),\overline v(x,t),\underline u(x,t)$ and $\underline v(x,t)$ in Definition \ref{def:def1} are continuous on $\mathbb{R}\times(-\infty,T_0)$, but not smooth at some points, then $(u_n(x,t),v_n(x,t))$ is still a unique classical bounded solution of \eqref{eq:gen} with $(u_n(x,-n),v_n(x,-n))=(\underline u(x,-n),\underline v(x,-n))$ for $x\in\mathbb{R}$ and $t>-n$, and \eqref{eq:gen} also admits an entire solution. Thus, in the proof of the existence of entire solutions of \eqref{eq:gen}, the smoothness of $\overline u,\overline v,\underline u$ and $\underline v$ in Definition \ref{def:def1} can be weakened into almost everywhere $C^1$ in $t$ and $C^2$ in $x$.
\end{remark}

\section{Existence of different types of entire solutions}
In this section, different methods of constructing pairs of coupled super-sub solutions are applied in order to obtain entire solutions. To begin with the discussion, from Lemma \ref{lem:asym1} to Lemma \ref{lem:asym05}, one can easily get the following lemma.
\begin{lemma}\label{lem:esti}
Suppose $(\phi(\xi),\psi(\xi))$ is the solution to \eqref{eq:OLV}-\eqref{eq:MC} under the case $(iv)$. When $\xi\leqslant0$, there exist some positive numbers $M_1$, $\overline{M}_1$, $M_2$, $\overline{M}_2$ and $\kappa$ such that
$$\overline{M}_1\leqslant\frac{\phi(\xi)}{\phi'(\xi)}\leqslant M_1,\ \ \ \overline{M}_2e^{\kappa\xi}\leqslant\phi(\xi)\leqslant M_2e^{\kappa\xi},\ \ \ \phi'(\xi)\leqslant M_2e^{\kappa\xi},$$
$$\overline{M}_1\leqslant\frac{\psi(\xi)}{\psi'(\xi)}\leqslant M_1,\ \ \ \overline{M}_2e^{\kappa\xi}\leqslant\psi(\xi)\leqslant M_2e^{\kappa\xi},\ \ \ \psi'(\xi)\leqslant M_2e^{\kappa\xi}.$$
On the other hand,  there are some positive constant $M_3$, $\overline{M}_3$, $M_4$, such that, for $\xi\geqslant0$,
$$\overline{M}_3\leqslant\frac{\phi'(\xi)}{\max\{u^{\ast}-\phi(\xi),v^{\ast}-\psi(\xi)\}}
\leqslant M_3,\ \ \ \ \phi'(\xi)\leqslant M_4e^{\lambda_2\xi},
$$ $$\overline{M}_3\leqslant\frac{\psi'(\xi)}{\max\{u^{\ast}-\phi(\xi),v^{\ast}-\psi(\xi)\}}
\leqslant M_3, \ \ \ \ \psi'(\xi)\leqslant M_4e^{\lambda_2\xi},$$ where $\lambda_2$ is given in Lemma \ref{lem:asym1}.
\end{lemma}

In addition, we define
$$\begin{array}{ll}
F_3(u,v)=\partial_tu-\partial_{xx}u-u(1-u-k_1v),\\[0.2cm]
F_4(u,v)=\partial_tv-d\,\partial_{xx}v-rv(1-v-k_2u).
\end{array}$$

\subsection{The first kind of entire solutions}
In this subsection, when solutions of \eqref{eq:LV} without diffusion are taken into account, different types of entire solutions to \eqref{eq:LV} are exhibited. Firstly we introduce some properties of the solution of the diffusion-free system
\begin{equation}\label{eq:odelv}
\left\{\begin{array}{lllll}
p'_1(t)=p_1(t)(1-p_1(t)-k_1q_1(t)),\\[0.2cm]
q'_1(t)=rq_1(t)(1-q_1(t)-k_2p_1(t)),\\[0.2cm]
\lim\limits_{t\rightarrow-\infty}(p_1(t),q_1(t))=(0,0),\\[0.2cm] \lim\limits_{t\rightarrow+\infty}(p_1(t),q_1(t))=(u^{\ast},v^{\ast}),\\[0.2cm]
(p_1(0),q_1(0))=(\theta_1,\theta_2)\in(0,u^{\ast})\times(0,v^{\ast}).
\end{array}
\right.
\end{equation}
It is not hard to see that $$(0,0)<\left(\frac{u^{\ast}\hat{\beta}_1e^{u^{\ast}t}}{1+\hat{\beta}_1e^{u^{\ast}t}},
\frac{v^{\ast}\hat{\beta}_2e^{v^{\ast}t}}{1+\hat{\beta}_2e^{v^{\ast}t}}\right)
\leqslant(p_1(t),q_1(t))\leqslant(u^{\ast},v^{\ast}),$$ where $\hat{\beta}_1=\frac{\theta_1}{u^{\ast}-\theta_1}$, $\hat{\beta}_2=\frac{\theta_2}{u^{\ast}-\theta_2}$, due to the Cauchy-Lipschitz theorem and $p'_1(t)>0$, $q'_1(t)>0$.

Then we will use the solutions of \eqref{eq:OLV} with \eqref{eq:BC} and \eqref{eq:MC}, namely the traveling front solutions of \eqref{eq:LV} connecting the origin and the positive equilibrium and their reflects as well as the solutions of \eqref{eq:odelv} to construct pairs of coupled super-sub solutions. Inspired by the papers \cite{mt09} and \cite{wl13}, we initially want to find the entire solutions similar with the annihilation type and use them to illustrate one species invading from both sides of the $x-$axis and coexisting with the other species. However, whether taking the method in \cite{wl13} or the way from \cite{mt09}, we can not get any pairs of coupled super-sub solutions of \eqref{eq:LV} satisfying \eqref{condi:csl} in the case (iv). As a result, we turn to find entire solutions similar to that in \cite{gm05} and \cite{hn99}. The way to construct super-sub solutions is similar to the one in \cite{wl15}, and different kinds of entire solutions have been obtained. One of them has the following asymptotic behavior. One component of this entire solution starts from 0 at $t\approx -\infty$ and as $t$ goes forward, the entire solution will stay in a conformed region. This thing implies one species invades from both sides of the $x-$axis and will finally mix with the other species.

It is not hard to prove the following pairs of functions
\begin{equation}\label{eq:sup-low1}
\begin{array}{llll}
\overline{u}(x,t)=1,
\\[0.2cm] \underline{u}_{ijm}(x,t)=\max\{\chi_i\phi(x+ct),\chi_j\phi(-x+ct),\chi_mp_1(t)\},
\\[0.2cm]
\overline{v}_{ijm}(x,t)=\min\{\chi_i\psi(x+ct),\chi_j\psi(-x+ct),\chi_mq_1(t)\},
\\[0.2cm]
\underline{v}(x,t)=0
\end{array}
\end{equation}
are pairs of coupled super-sub solutions of \eqref{eq:LV} defined on $\mathbb{R}\times\mathbb{R}$, where $\chi_i=i~(i=0,1)$, $\chi_j=j~(j=0,1)$, $\chi_m=m~(m=0,1)$, and $i$, $j$, $m$ are not zero at the same time.
\begin{theorem}\label{thm:entire1}
Under the case (iv), suppose that $(\phi(\xi),\psi(\xi))$ are the solution of \eqref{eq:OLV} with \eqref{eq:BC} and \eqref{eq:MC}, where $c\geqslant2\max\{1,\sqrt{rd}\}$, and $(p_1(t),q_1(t))$ are the solution of \eqref{eq:odelv}. Then \eqref{eq:LV} admits entire solutions $(u_{ijm}(x,t),v_{ijm}(x,t))$ on $\mathbb{R}\times\mathbb{R}$, satisfying
$(\underline{u}_{ijm},\underline{v})\leqslant(u_{ijm},v_{ijm})
\leqslant(\overline{u},\overline{v}_{ijm})$, where $(\underline{u}_{ijm},\underline{v}), (\overline{u},\overline{v}_{ijm})$ are given in (\ref{eq:sup-low1}), and the following properties:\\
${\rm(i)}$ $(u_{110}(x,t),v_{110}(x,t))=(u_{110}(-x,t),v_{110}(-x,t))$;\\
${\rm(ii)}$ $
\lim\limits_{t\rightarrow-\infty}\sup\limits_{x\in\mathbb{R}}|v_{110}(x,t)|=0;
$\\
${\rm(iii)}$ for any bounded and closed intervals $\mathbb{I}$ of $\mathbb{R}$,
$$\lim\limits_{x\rightarrow\pm\infty}
\sup\limits_{t\in\mathbb{I}}|v_{110}(x,t)|=0;$$\\
${\rm(iv)}$ $$ u^{\ast}\leqslant\lim\limits_{t\rightarrow+\infty}\sup\limits_{x\in\mathbb{R}}
|u_{110}(x,t)|
\leqslant1,\qquad 0\leqslant\lim\limits_{t\rightarrow+\infty}\sup\limits_{x\in\mathbb{R}}
|v_{110}(x,t)|
\leqslant v^{\ast}.$$
\end{theorem}

\Proof The existence of entire solutions can be obtained from Theorem \ref{thm:exig}, thus here we focus on the properties of entire solutions.

From \eqref{eq:sup-low1} and Theorem \ref{thm:exig}, the properties (i) and (ii) are obvious. The proofs of the properties (iv) and (v) are simple and we omit it here.
\qed

\begin{remark}\label{remark44}
In \cite{mt09} and \cite{wl15}, the authors used some exact solutions of \eqref{eq:LV} to construct the super-sub solutions to obtain the existence of entire solution. Although, under the case (iv), we also can find the exact solutions of \eqref{eq:LV} connecting the origin and the positive equilibrium in the paper \cite{h12}, we can not get any pairs of coupled super-sub solutions followed the method in \cite{mt09} or \cite{wl15}. Thus we here use the traveling front solutions and their reflects to construct the super-sub solutions to establish entire solutions. When $(i,j,m)=(1,1,0)$, the entire solution describes that the species invades from both sides of the $x-$axis and finally mix with the other species.
\end{remark}

\subsection{The second kind of entire solutions}
Finally we use the solutions from different scalar equations to construct the super-sub solutions of \eqref{eq:LV} which leads to the existence of entire solutions.

When $0<k_1,k_2<1$, from the paper \cite{kpp37} it is not hard to prove that the equation
$$\partial_tu-\partial_{xx}u-u(1-k_1-u)=0$$
admits traveling front solutions $u=\phi(\xi_1)$, $\xi_1=x+s_1t$, $s_1\geqslant2\sqrt{1-k_1}$, connecting 0 and $1-k_1$, and the equation
$$\partial_tv-d\,\partial_{xx}v-rv(1-k_2-v)=0$$
also admits traveling front solutions $v=\psi(\xi_2)$, $\xi_2=x+s_2t$, $s_2\geqslant2\sqrt{r(1-k_2)}$, connecting 0 and $1-k_2$. That is to say, $\phi(\xi_1)$ and $\psi(\xi_2)$ respectively satisfy
\begin{equation}\label{eq:sode1}
\phi''-c\phi'+\phi(1-k_1-\phi)=0,\quad\lim\limits_{\xi_1
\rightarrow-\infty}\phi(\xi_1)=0,\quad\lim\limits_{\xi_1
\rightarrow+\infty}\phi(\xi_1)=1-k_1,
\end{equation}
and
\begin{equation}\label{eq:sode2}
d\psi''-c\psi'+r\psi(1-k_2-\psi)=0,\ \ \lim\limits_{\xi_2
\rightarrow-\infty}\psi(\xi_2)=0,\ \ \lim\limits_{\xi_2
\rightarrow+\infty}\psi(\xi_2)=1-k_2.
\end{equation}
Furthermore, from the paper \cite{fm77}, $\phi'(\xi_1),\psi'(\xi_2)>0$. Then we will prove the following theorem.
\begin{theorem}\label{thm:entire2}
Suppose that $0<k_1,k_2<1$, and $\phi(\xi_1)$, $\psi(\xi_2)$ satisfies \eqref{eq:sode1}
and \eqref{eq:sode2} respectively, where $s_1\geqslant2\sqrt{1-k_1}$, $s_2\geqslant2\sqrt{r(1-k_2)}$. Then \eqref{eq:LV} admits an entire solution $(u_1(x,t),v_1(x,t))$ on $\mathbb{R}\times\mathbb{R}$ satisfying
$(\underline{u},\underline{v})\leqslant(u_1,v_1)
\leqslant(\overline{u},\overline{v})$,
where
\begin{align*}
\overline{u}(x,t)=1,\qquad\underline{u}(x,t)=
\max\{\phi(x+s_1t),\phi(-x+s_1t)\},\\[0.2cm]
\overline{v}(x,t)=1,\qquad\underline{v}(x,t)=
\max\{\psi(x+s_2t),\psi(-x+s_2t)\}.
\end{align*}
In addition, on $\mathbb{R}\times\mathbb{R}$, $u_1(x,t)=u_1(-x,t)$, $v_1(x,t)=v_1(-x,t)$.
\end{theorem}
\Proof By subsitituting $(\overline{u},\underline{v})$ into $F_3(u,v)$ and $F_4(u,v)$, it is not hard to prove
$F_3(\overline{u},\underline{v})=k_1\underline{v}\geqslant0$,
$F_4(\overline{u},\underline{v})=c\underline{v}'-d\underline{v}''-r
\underline{v}(1-k_2-\underline{v})\leqslant0$.
Similarly, $F_3(\underline{u},\overline{v})=c\underline{u}'-\underline{u}''-
\underline{u}(1-k_1-\underline{u})\leqslant0$,
$F_4(\underline{u},\overline{v})=
k_2\underline{u}\geqslant0$. The rest of the proof can be derived from Theorem \ref{thm:exig}.\qed


\section*{References}
\bibliographystyle{elsarticle-num}

\begin{thebibliography}{aa}
\bibitem{al91} S. Ahmad, A. C. Lazer, An elementary approach to traveling front solutions to a system of N competition-diffusion equations, Nonlinear Anal. 16 (1991) 892-901.

\bibitem{alt08} S. Ahmad, A. C. Lazer, A. Tineo, Traveling waves for a system of equations, Nonlinear Anal. 68 (2008) 3909-3912.

\bibitem{cg05} X. Chen, J. S. Guo, Existence and uniqueness of entire solutions for a reaction-diffusion equation, J. Differential Equations 212 (2005) 62-84.

\bibitem{cg84} C. Conley, R. Gardner,  Application of the generalized Morse index to travelling wave solutions of a competitive reaction-diffusion model, Indiana Univ. Math. J. 33 (1984) 319-343.

\bibitem{ct12} E.C.M. Crooks, J. C. Tsai, Front-like entire solutions for equations with convection, J. Differential Equations 253 (2012) 1206-1249.

\bibitem{fm77} P. C. Fife, J. B. McLeod, The approach of solutions of nonlinear diffusion equations to travelling front solutions, Arch. Ration. Mech. Anal., {\bf 65} (1977) 335-361.

\bibitem{fmn04} Y. Fukao, Y. Morita, H. Ninomiya, Some entire solutions of the Allen-Cahn equation, Taiwanese J. Math. 8 (2004) 15-32.

\bibitem{gm05} J. S. Guo, Y. Morita, Entire solutions of reaction-diffusion equations and an application to discrete diffusive equations, Discrete Contin. Dyn. Syst. 12 (2005) 193-212.

\bibitem{hn99} F. Hamel, N. Nadirashvili, Entire solutions of the KPP Equation,
Comm. Pure Appl. Math. 52 (1999) 1255-1276.

\bibitem{hn01} F. Hamel, N. Nadirashvili, Travelling fronts and entire solutions
of the Fisher-KPP Equation in R$^N$, Arch. Rational Mech. Anal. 157 (2001) 91-163.

\bibitem{h12} L. C. Hung, Exact traveling wave solutions for diffusive Lotka-Volterra systems of two competing species, Japan J. Indust. Appl. Math., {\bf 29} (2012) 237-251.

\bibitem{kpp37} A. Kolmogorov, I. Petrovsky and N. Piscounoff, \'{E}tude de \'{I}equation de la diffusion avec croissance de la quantit\'{e} de mati\`{e}re et son application \`{a} unprobleme biologique, Bull. Univ. Moskou, Ser. Internat., Sec. A, {\bf 1} (1937) 6, 1-25.
    
\bibitem{lhf11} Anthony W. Leung, X. J. Hou, W. Feng, Traveling wave solutions for Lotka-Volterra system re-visited, Discret. Contin. Dyn. Syst. Ser. B 15 (2011) 171-196.

\bibitem{lhl08} Anthony W. Leung, X. J. Hou, Y. Li, Exclusive traveling waves for competitive reaction-diffusion systems and their stabilities, J. Math. Anal. Appl. 338 (2008) 902-924.

\bibitem{ll12} K. Li, X. Li, Asymptotic behavior and uniqueness of traveling wave solutions in Ricker competition system, J. Math. Anal. Appl. 389 (2012) 486-497.

\bibitem{llw08} W. T. Li, N. W. Liu, Z. C. Wang, Entire solutions in reaction-advection-diffusion equations in cylinders, J. Math. Pures Appl. 90 (2008) 492-504.

\bibitem{lww08} W. T. Li, Z. C. Wang, J. H. Wu, Entire solutions in monostable reaction-diffusion equations with delayed nonlinearity, J. Differential Equations 245 (2008) 102-129.

\bibitem{lzz15} W. T. Li, L. Zhang, G. B. Zhang, Invasion entire solutions in a competition system with nonlocal dispersal, Discret. Contin. Dyn. Syst. 35 (2015) 1531-1560.

\bibitem{l08} Z. Y. Li, Asymptotic behavior of traveling wave fronts of Lotka-Volterra competitive system, Int. Journal of Math. Analysis 2 (2008) 1295-1300.

\bibitem{l12} G. Y. Lv, Entire solutions of delayed reaction diffusion equations, Z. Angew. Math. Mech. 92 (2012) 204-216.

\bibitem{m79} P de Mottoni, Qualitative analysis for some quasilinear parabolic systems, Inst. Math. Polish Acad. Sci. Zam 190 (1979) 11-79.

\bibitem{mn06} Y. Morita, H. Ninomiya, Entire solutions with merging fronts to reaction-diffusion equations, J. Dynam. Differential Equations 18 (2006) 841-861.

\bibitem{mt09} Y. Morita, K. Tachibana, An entire solution to the Lotka-Volterra competition-diffusion equations, SIAM J. Math. Anal. 40 (2009) 2217-2240.

\bibitem{tf80} M. M. Tang, P. C. Fife, Propagating fronts for competing species equations with diffusion, Arch. Rational Mech. Anal. 73 (1980) 69-77.

\bibitem{t82} D. Terman, Comparison theorems for reaction-diffusion systems defined in an unbounded domain, Technical Summary Report, 1982.

\bibitem{v95} J. H. V. Vuuren, The existence of travelling plane waves in a general class of competition-diffusion systems, IMA J. Appl. Math. 55 (1995) 135-148.

\bibitem{wl10} M. X. Wang, G. Y. Lv, Entire solutions of a diffusive and competitive Lotka-Volterra type system with nonlocal delays, Nonlinearity 23 (2010) 1609-1630.

\bibitem{wl15} Y. Wang, X. Li, Some entire solutions to the competitive reaction diffusion system, J. Math. Anal. Appl. 430 (2015) 993-1008.

\bibitem{wl13} X. H. Wang, G. Y. Lv, Entire solutions for Lotka-Volterra competition-diffusion model, Int. J. Biomath. 6 (2013) 1350020.

\bibitem{wlr09} Z. C. Wang, W. T. Li, S. Ruan, Entire solutions in bistable reaction-diffusion equations with nonlocal delayed nonlinearity, Trans. Am. Math. Soc. 361 (2009) 2047-2084.

\bibitem{wlw09} Z. C. Wang, W. T. Li, J. Wu, Entire solutions in delayed lattice differential equations with monostable nonlinearity, SIAM J. Math. Anal. 40 (2009) 2392-2420.

\bibitem{ww13} S. L. Wu, H. Y. Wang, Front-like entire solutions for monostable reaction-diffusion systems, J. Dynam. Differential Equations 25 (2013) 505-533.

\bibitem{y03} H. Yagisita, Backward global solutions characterizing annihilation dynamics of travelling fronts, Publ. Res. Inst. Math. Sci 39 (2003) 117-164.

\bibitem{yy11} Z. X. Yu, R. Yuan, Traveling waves for a Lotka-Volterra competition system with diffusion, Mathematical and Computer Modelling 53 (2011) 1035-1043.

\end{thebibliography}

\end{document}